\documentclass[12pt]{amsart}
\usepackage{fullpage,url}
\usepackage[dvips]{hyperref}
 \usepackage{amscd}   
\usepackage{amsmath}

\DeclareFontEncoding{OT2}{}{} 

\def\dsp{\def\baselinestretch{2.0}\large\normalsize}
\dsp

\newcommand{\F}{\mathbb{F}}

\newcommand{\Q}{\mathbb{Q}}

\DeclareMathOperator{\Span}{Span}
\DeclareMathOperator{\per}{per}
\DeclareMathOperator{\alg}{alg}

\DeclareMathOperator{\ybar}{\overline{y}}
\DeclareMathOperator{\xbar}{\overline{x}}

\DeclareMathOperator{\ac}{ac}
\DeclareMathOperator{\tor}{tor}
\DeclareMathOperator{\hhat}{\widehat{h}}
\DeclareMathOperator{\Hhat}{\widehat{H}}
\DeclareMathOperator{\trdeg}{tr deg}

\DeclareMathOperator{\ord}{ord}

\DeclareMathOperator{\Frac}{Frac}

\newcommand{\tensor}{\otimes}

\newtheorem{theorem}{Theorem}[section]
\newtheorem{lemma}[theorem]{Lemma}
\newtheorem{corollary}[theorem]{Corollary}

\newtheorem{Fact}[theorem]{Fact}

\theoremstyle{definition}
\newtheorem{definition}[theorem]{Definition}

\newtheorem{conjecture}[theorem]{Conjecture}
\newtheorem{example}[theorem]{Example}
\newtheorem{Claim}[theorem]{Claim}

\theoremstyle{remark}
\newtheorem{remark}[theorem]{Remark}

\newtheorem{sublemma}[theorem]{Sublemma}

\title{The Lehmer inequality and the Mordell-Weil theorem for Drinfeld modules}
\author{Dragos Ghioca}

\begin{document}
\maketitle
\address{Department of Mathematics and Statistics, Hamilton Hall, Room 218, McMaster University, 1280 Main Street West, Hamilton, Ontario L8S 4K1, Canada \\ \url{dghioca@math.mcmaster.ca}\\ \emph{phone}: 905-525-9140, EXT. 26079\\ \emph{fax}: 905-522-0935}
\newpage
\begin{abstract}
In this paper we prove several Lehmer type inequalities for Drinfeld 
modules which will enable us to prove certain Mordell-Weil type structure theorems for Drinfeld modules.

\end{abstract}

\newpage
\section{Introduction}
The classical Lehmer conjecture (see \cite{Leh}, page $476$) asserts that there is
an absolute constant $C>0$ 
so that any algebraic number $\alpha$ that is not a
root of unity satisfies the following inequality for its logarithmic height
$$h(\alpha)\ge\frac{C}{[\Q(\alpha):\Q]}.$$
A partial result 
towards this conjecture is obtained in \cite{Dob}. The
analog of Lehmer 
conjecture for elliptic curves and abelian varieties asks for a good lower bound for the canonical height of a non-torsion point of the abelian variety. Also this question has been much studied 
(see \cite{Siabba},
\cite{DaHi}, \cite{HiSi}, \cite{Mass}, \cite{Siabell}). 

In this paper we prove several inequalities for the height of non-torsion points of Drinfeld modules. These inequalities have the same flavor as the above mentioned Lehmer's conjecture. Using our inequalities we will be able to prove several Mordell-Weil type structure theorems for Drinfeld modules over certain infinitely generated fields. Next we will define the notion of Drinfeld modules.

In this paper we will use the following notation: $p$ is a prime number and $q$ is 
a power of $p$. We denote by $\F_q$ the finite field with $q$ elements. 
We let $C$ be a nonsingular projective curve defined over $\F_q$ and we 
fix a closed point $\infty$ on $C$. Then we define $A$ as the ring of functions on $C$ that are regular 
everywhere except possibly at $\infty$.

We let $K$ be a field extension of $\F_q$. We fix a 
morphism $i:A\rightarrow K$. We define the operator $\tau$ as the power of the 
usual Frobenius with the property that for every $x\in K^{\alg}$, $\tau(x)=x^q$. 
Then we let $K\{\tau\}$ be the ring of polynomials in $\tau$ with coefficients from $K$ (the addition is the usual one while the multiplication is the composition of functions).

We fix an algebraic closure of $K$, denoted $K^{\alg}$. We denote by $\F_p^{\alg}$ the algebraic closure of $\F_p$ inside $K^{\alg}$. Also, for us, the symbol "$\subset$" means inclusion, not neccessarily strict inclusion. 

A Drinfeld module is a morphism $\phi:A\rightarrow K\{\tau\}$ for which the 
coefficient of $\phi_a$ is $i(a)$ for every $a\in A$, and there exists $a\in A$ such that $\phi_a\ne i(a)\tau^0$. Following the definition 
from \cite{Goss} we call $\phi$ a Drinfeld module of generic characteristic 
if $\ker(i)=\{0\}$ and we call $\phi$ a Drinfeld module of finite 
characteristic if $\ker(i)\ne \{0\}$.

For each field $L$ containing $K$, $\phi(L)$ denotes $L$ with the $A$-action 
given by $\phi$. When $K$ is a finitely generated extension of 
$\F_q$, it was proved in \cite{Poo} (in the case that $\trdeg_{\F_p}K=1$) 
and in \cite{Wan} (for arbitrary finite, positive transcendence degree) that $\phi(K)$ 
is the direct sum of a finite torsion submodule with a free submodule of rank 
$\aleph_0$. We will prove in Theorem ~\ref{T:T4mw} a similar structure result for certain infinitely generated extensions of $\F_q$.

The key result is Theorem \ref{T:T2mwg}, which will be proved in the fourth
section of this paper. Mainly what will be needed for our result will be a 
better understanding of the heights associated to $\phi$, both local and global 
heights. These heights were first introduced in \cite{Den} and then 
contributions towards their understanding were done in \cite{teza}, \cite{Poo} and \cite{Wan}. We mention that our results are part of our Ph.D. thesis \cite{teza}. Lemma~\ref{L:L5} appears also in our paper \cite{Cre}, in which we prove a local Lehmer inequality for Drinfeld modules.

Theorem \ref{T:T2mwg} will also give us the technical ingredient to obtain an 
uniform boundedness result for the torsion submodule of $\phi(K)$. This will be explained in Corollary ~\ref{C:CorT2}.

\section{Heights associated to a Drinfeld module}

We continue with the notation from Section $1$. So, $K$ is a field extension of $\F_q$ and $\phi:A\rightarrow K\{\tau\}$ is a Drinfeld module. We define $M_K$ as the set of all discrete valuations 
of $K$. We also normalize all the valuations $v\in M_K$ such that the range of $v$ is $\mathbb{Z}$. 

\begin{definition}
\label{D:good sets}
We call a subset $U\subset M_K$ equipped with a function 
$d:U\rightarrow\mathbb{R}_{>0}$ a \emph{good set of valuations}
if the following properties are satisfied

(i) for every nonzero $x\in K$, there are finitely
many $v\in U$ such that $v(x)\ne 0$.

(ii) for every nonzero $x\in K$,
$$\sum_{v\in U}d(v)\cdot v(x)=0.$$
The positive real number $d(v)$ will be called the \emph{degree} of the 
valuation $v$. When we say that the positive real number $d(v)$ is associated to 
the valuation $v$, we understand that the degree of $v$ is $d(v)$. 

When $U$ is a good set of valuations, we will refer to
property (ii) as the sum formula for $U$.
\end{definition}

\begin{definition}
\label{D:constants}
Let $U$ be a good set of valuations on $K$. The set
$\{0\}\cup\{x\in K\mid v(x)=0\text{ for all }v\in U\}$
is the set of \emph{constants} for $U$. We denote this
set by $C(U)$.
\end{definition}

\begin{lemma}
\label{L:Lconstants}
Let $U$ be a good set of valuations on $K$. If $x\in
K$ is integral at all places $v\in U$, then $x\in
C(U)$.
\end{lemma}

\begin{proof}
Let $x\in K\setminus\{0\}$. By the sum formula for
$U$, if $v(x)\ge 0$ for all $v\in U$, then actually
$v(x)=0$ for all $v\in U$ (a sum of non-negative
numbers is $0$ if and only if all the numbers are
$0$).
\end{proof}

\begin{lemma}
\label{L:field of constants}
Let $U$ be a good set of valuations on a field $K$.
The set $C(U)$ is a subfield of $K$.
\end{lemma}

\begin{proof}
By its definition, $C(U)$ is closed under
multiplication and division. Because of Lemma
~\ref{L:Lconstants}, $C(U)$ is closed also under
addition.
\end{proof}

\begin{definition}
\label{D:loccoh}
Let $v\in M_K$ of degree $d(v)$. We say that the valuation $v$ is coherent (on 
$K^{\alg}$) if for every finite extension $L$ of $K$,
\begin{equation}
\label{E:(iii)}
\sum_{\substack{w\in M_{L}\\
w|v}}e(w|v)f(w|v)=[L:K],
\end{equation}
where $e(w|v)$ is the ramification index and $f(w|v)$ is the relative degree 
between the residue field of $w$ and the residue field of $v$.

Condition \eqref{E:(iii)} says that $v$ is \emph{defectless} in $L$. In this 
case, we also let the degree of any $w\in M_L$, $w|v$ be
\begin{equation}
\label{E:(ii)}
d(w)=\frac{f(w|v)d(v)}{[L:K]}.
\end{equation}
\end{definition}

It is immediate to see that \eqref{E:(ii)} is equivalent
to the stronger condition that for every two finite
extensions of $K$, $L_1\subset L_2$ and for every
$v_2\in M_{L_2}$ that lies over $v_1\in M_{L_1}$, which in turn lies over $v$, 
\begin{equation}
\label{E:towernew}
d(v_2)=\frac{f(v_2|v_1)d(v_1)}{[L_2:L_1]}.
\end{equation}

We will use in our proofs the following result from \cite{End} (see $(18.1)$, 
page $136$).
\begin{lemma}
\label{L:important valuation lemma}
Let $L_1\subset L_2\subset L_3$ be a tower of finite
extensions. Let $v\in M_{L_1}$ and denote by
$w_1,\dots,w_s$ all the places of $L_2$ that lie over
$v$. Then the following two statements are equivalent:

1) $v$ is defectless in $L_3$.

2) $v$ is defectless in $L_2$ and $w_1,\dots,w_s$ are
defectless in $L_3$.
\end{lemma}

Lemma \ref{L:important valuation lemma} shows that
condition \eqref{E:(iii)} of Definition ~\ref{D:loccoh} is equivalent to the 
following statement: for
every two finite extensions of $K$, $L_1\subset L_2$
and for every $v_1\in M_{L_1}$, $v_1|v$
\begin{equation}
\label{E:L_2L_1}
\sum_{\substack{w\in M_{L_2}\\
w|v_1}}e(w|v_1)f(w|v_1)=[L_2:L_1].
\end{equation}
The following result is an immediate consequence of Definition \ref{D:loccoh} 
and Lemma~\ref{L:important valuation lemma}.
\begin{lemma}
\label{L:always coherent}
If $v\in M_K$ is a coherent valuation (on $K^{\alg}$), then for every finite 
extension $L$ of $K$ and for every $w\in M_L$ and $w|v$, $w$ is a coherent 
valuation (on $K^{\alg}=L^{\alg}$).
\end{lemma}

\begin{definition}
\label{D:coherent good sets}
We let $U_K$ be a good set of valuations on $K$. We call $U_K$ a \emph{coherent} 
good set of
valuations (on $K^{\alg}$) if the following two conditions are satisfied

(i) for every finite extension $L$ of $K$, if $U_L\subset M_L$ is the set of all 
valuations lying over valuations from $U_K$, then $U_L$ is a good set of 
valuations.

(ii) for every $v\in U_K$, the valuation $v$ is coherent (on $K^{\alg}$).
\end{definition}

\begin{remark}
\label{R:important remark}
Using the argument from page $9$ of \cite{Ser}, we conclude that condition (i) 
from Definition ~\ref{D:coherent good sets} is automatially satisfied if $U_K$ 
is a good set of valuations and if condition (ii) of Definition ~\ref{D:coherent 
good sets} is satisfied.
\end{remark} 

An immediate corollary to Lemma \ref{L:always coherent} is the following result.
\begin{corollary}
\label{C:always coherent}
If $U_K\subset M_K$ is a good set of valuations that is coherent (on 
$K^{\alg}$), then for every finite extension $L$ of $K$, if $U_L$ is the set of 
all valuations on $L$ which lie over valuations from $U_K$, then $U_L$ is a 
coherent good set of valuations. 
\end{corollary}

Fix now a field $K$ of characteristic $p$ and let
$\phi:A\rightarrow K\{\tau\}$ be a Drinfeld module.
Let $v\in M_K$ be a coherent valuation (on $K^{\alg}$). Let $d(v)$ be the degree 
of $v$ as in Definition ~\ref{D:loccoh}. For such $v$, we construct the local 
height $\hhat_{v}$ with respect to the Drinfeld 
module $\phi$. Our construction follows \cite{Poo}.
For $x\in K$, we set $\tilde{v}(x)=\min\{
0,v(x)\}$. For a non-constant element $a\in A$, we
define 
\begin{equation}
\label{E:localheight1}
V_v(x)=\lim_{n\rightarrow\infty}\frac{\tilde{v}(\phi_{a^n}(x))}{\deg(\phi_{a^n})
}.
\end{equation}
This function is well-defined and satisfies the same properties as in
Propositions $1$-$3$ from 
\cite{Poo}. Mainly, we will use the following facts:

$1)$ if $x$ and all the coefficients of $\phi_a$ are integral at $v$, then 
$V_v(x)=0$.

$2)$ for all $b\in A\setminus\{0\}$,
$V_v(\phi_b(x))=\deg(\phi_b)\cdot V_v(x)$. Moreover,
we can use any non-constant $a\in A$ for the
definition of $V_v(x)$ and we will always get the same
function $V_v$.

$3)$ $V_v(x\pm y)\ge\min\{V_v(x),V_v(y)\}$.

$4)$ if $x\in\phi_{\tor}$, then $V_v(x)=0$.

We define then
\begin{equation}
\label{E:localheight2}
\hhat_{v}(x)=-d(v) V_v(x).
\end{equation}

If $L$ is a finite extension of $K$ and $w\in M_L$ lies over $v$ then we define 
similarly the function $V_w$ on $L$ and just as above, we let 
$\hhat_w(x)=-d(w)V_w(x)$ for every $x\in L$.

If $U=U_K\subset M_K$ is a coherent good set of valuations, then for each $v\in 
U$, 
we denote by $\hhat_{U,v}$ the local height associated to $\phi$ with respect to 
$v$ (the construction of $\hhat_{U,v}$ is identical with the one from above). 
Then we define the global height associated to $\phi$ as
\begin{equation}
\label{E:globalheight}
\hhat_U(x)=\sum_{v\in U}\hhat_{U,v}(x).
\end{equation}
For each $x$, the above sum is finite due to fact $1)$ stated above
(see also Proposition $6$ of \cite{Poo}). 

For each finite extension $L$ of $K$, we let $U_L$ be the set of all valuations 
of $L$ that lie over places from $U_K$. As stated in Corollary \ref{C:always 
coherent}, $U_L$ is also a coherent good set of valuations and so, we can define 
the local heights $\hhat_{U_L,w}$ with respect to $w\in U_L$, associated to $\phi$ for all 
elements $x\in L$. Then we define the global height of $x$ as
$$\hhat_{U_L}(x)=\sum_{w\in U_L}\hhat_{U_L,w}(x).$$

\begin{Claim}
\label{C:coerent}
Let $L_1\subset L_2$ be finite extensions of $K$. Let $v\in U_{L_1}$ and $x\in 
L_1$. Then 
$$\sum_{\substack{w\in 
U_{L_2}\\w|v}}\hhat_{U_{L_2},w}(x)=\hhat_{U_{L_1},v}(x).$$
\end{Claim}

\begin{proof}
We have 
$$\sum_{\substack{w\in U_{L_2}\\w|v}}\hhat_{U_{L_2},w}(x)=-\sum_{\substack{w\in 
U_{L_2}\\w|v}}d(w)V_w(x).$$
Because $d(w)=\frac{d(v)f(w|v)}{[L_2:L_1]}$ (see \eqref{E:towernew}) and 
$V_w(x)=e(w|v)V_v(x)$ we get
$$\sum_{\substack{w\in 
U_{L_2}\\w|v}}\hhat_{U_{L_2},w}(x)=\frac{-d(v)V_v(x)}{[L_2:L_1]}\sum_{\substack{
w\in U_{L_2}\\w|v}}e(w|v)f(w|v).$$
Because $v$ is defectless and $\hhat_{U_{L_1},v}(x)=-d(v)V_v(x)$, we are done.
\end{proof}

Claim \ref{C:coerent} shows that our 
definition of the global height is independent of the
field $L$ containing $x$ and so, we can drop the index
referring to the field $L$ containing $x$ when we work
with the global height associated to a coherent good
set of valuations.

The above construction for global heights
depends on the selected good set of valuations $U_K$
on $K$. If we work with global heights only for points $x\in K$, then $U_K$ can be any good set of valuations on $K$. If we are interested in global heights for all points $x\in K^{\alg}$, then $U_K$ has to be a good set of valuations on $K$, which is coherent (on $K^{\alg}$). Also, technically speaking, we do not need for local heights $\hhat_v$ the valuation $v$ be coherent as long as we restrict ourselves to points $x\in K$. We will always specify first which is the good
set of valuations that we consider when we will work
with heights associated to a Drinfeld module.

\section{Examples of coherent good sets of valuations}
Let $F$ be a field of characteristic $p$ and let
$K=F(x_1,\dots,x_n)$ be the rational function field of
transcendence degree $n\ge 1$ over $F$. We let
$F^{\alg}$ be the algebraic closure of $F$ inside
$K^{\alg}$. We will construct a coherent good set of
valuations on $K$. 

First we construct a good set of valuations on $K$ and then we will show that 
this set is also coherent. According to Remark ~\ref{R:important remark}, we 
only need to show that each of the valuations on $K$ we construct is coherent.

Let $R=F[x_1,\dots,x_n]$. For each irreducible
polynomial $P\in R$ we let $v_P$ be the valuation on
$K$ given by 
$$v_P(\frac{Q_1}{Q_2})=\ord_P(Q_1)-\ord_P(Q_2)\text{
for every nonzero $Q_1,Q_2\in R$,}$$
where by $\ord_P(Q)$ we denote the order of the
polynomial $Q\in R$ at $P$.

Also, we construct the valuation $v_{\infty}$ on $K$
given by 
$$v_{\infty}(\frac{Q_1}{Q_2})=\deg (Q_2)-\deg
(Q_1)\text{ for every nonzero $Q_1,Q_2\in R$,}$$
where by $\deg(Q)$ we denote the total degree of the
polynomial $Q\in R$.

We let $M_{K/F}$ be the set of all valuations $v_P$
for irreducible polynomials $P\in R$ plus the
valuation $v_{\infty}$. We let the degree of $v_P$ be $d(v_P)=\deg(P)$ for
every irreducible polynomial $P\in R$ and we also let
$d(v_{\infty})=1$. Then, for every nonzero $x\in K$,
we have the sum formula
$$\sum_{v\in M_{K/F}}d(v)\cdot v(x)=0.$$
So, $M_{K/F}$ is a good set of valuations on $K$
according to Definition ~\ref{D:good sets}. The field
$F$ is the field of constants with respect to
$M_{K/F}$.

\begin{remark}
\label{R:geometrie}
The valuations constructed above are 
exactly the valuations associated with the 
irreducible divisors of the projective space $\mathbb{P}_{F}^n$. The degrees of 
the valuations are the projective degrees of the corresponding irreducible 
divisors.
\end{remark}

Let $K'$ be a finite extension of $K$ and let $F'$ be
the algebraic closure of $F$ in $K'$. We let
$M_{K'/F'}$ be the set of all valuations on $K'$ that
extend the valuations from $M_{K/F}$. We normalize
each valuation $w$ from $M_{K'/F'}$ so that the range
of $w$ is $\mathbb{Z}$. Also, we define 
\begin{equation}
\label{E:coherent degree}
d(w)=\frac{f(w|v)d(v)}{[K':K]}
\end{equation}
for every $w\in M_{K'/F'}$ and $v\in M_{K/F}$ such
that $w|v$. Note that strictly speaking, $w$ is an extension of $v$ as a place 
and \emph{not} as a valuation function. However, we still call $w$ an extension 
of $v$.

\begin{remark}
\label{R:geometrie2}
Continuing the observations made in Remark \ref{R:geometrie}, the valuations 
defined on $K'$ are the ones associated with irreducible divisors of the 
normalization of $\mathbb{P}_{F}^n$ in $K'$. In general, the discrete valuations 
associated with the irreducible divisors of a variety which is regular in 
codimension $1$ form a coherent good set of valuations.
\end{remark}

In order to show that $M_{K/F}$
is a coherent good set of valuations (on $K^{\alg}$),
we need to check that condition \eqref{E:(iii)} of Definition
~\ref{D:loccoh} is satisfied. This is proved in Chapter $1$, Section $4$ of 
\cite{sere} (Hypothesis (F) holds for algebras of finite type over fields and 
so, it holds for localizations of such algebras). For each $v\in M_{K/F}$ we 
apply Propositions $10$ and $11$ of \cite{sere} to the local ring of $v$ to show 
$v$ is coherent.

Now, in general, let $F$ be a field of characteristic
$p$ and let $K$ be any finitely generated extension over
$F$, of positive transcendence degree over $F$. If $F$
is algebraically closed in $K$, we construct a coherent good
set of valuations $M_{K/F}\subset M_K$, as follows. We
pick a transcendence basis $\{x_1,\dots,x_n\}$ for
$K/F$ and first construct as before the set of
valuations on $F(x_1,\dots,x_n)$:
$$\{v_{\infty}\}\cup\{v_P\mid P\text{ irreducible
polynomial in }F[x_1,\dots,x_n]\}.$$
Then, by Corollary \ref{C:always coherent}, we have a unique way of extending 
coherently
this set of valuations to a good set of valuations on
$K$. The set $M_{K/F}$ depends on our initial choice
of the transcendence basis for $K/F$. Thus, in our
notation $M_{K/F}$, we suppose that $K/F$ comes
equipped with a choice of a transcendence basis for
$K/F$.

We also note that for every $v\in M_{K/F}$, if $v_0\in
M_{F(x_1,\dots,x_n)/F}$ lies below $v$, then
\begin{equation}
\label{E:important inequality}
d(v)=\frac{f(v|v_0)d(v_0)}{\left[K:F(x_1,\dots,x_n)\right]}\ge
\frac{1}{\left[K:F(x_1,\dots,x_n)\right]}.
\end{equation}
In general, if $K'$ is a finite extension of $K$ and $v'\in M_{K'}$ lies above 
$v\in M_K$, then
\begin{equation}
\label{E:importanta}
d(v')=\frac{f(v'|v)d(v)}{[K':K]}\ge\frac{d(v)}{[K':K]}.
\end{equation}

For each such good set of valuations $M_{K/F}$ and for
any Drinfeld module $\phi:A\rightarrow K\{\tau\}$, we
construct as before the set of local heights and the global
height associated to $\phi$. We denote the local
heights by $\hhat_{M_{K/F},v}$ and the global height
by $\hhat_{M_{K/F}}$. If $F$ is a finite field, our
construction coincides with the one from \cite{Wan}.
Thus, if $F$ is a finite field, we will drop the
subscript $M_{K/F}$ from the notation of the local
heights and of the global height. Also, when $F$ is a
finite field and $\trdeg_{F}K=1$, our construction
also coincides with the one from \cite{Poo}.

\section{Lehmer inequality for Drinfeld modules}

The paper \cite{DenL} formulated a conjecture whose general form is Conjecture ~\ref{C:Con}, which we refer to as the Lehmer inequality for Drinfeld modules.
\begin{conjecture}
\label{C:Con}
Let $K$ be a finitely generated field. For any Drinfeld module $\phi:A\rightarrow K\{\tau\}$ there exists a constant 
$C>0$ depending only on $\phi$ such that any non-torsion point $x\in K^{\alg}$ satisfies 
$\hhat(x)\ge\frac{C}{[K(x):K]}$.
\end{conjecture}

We fix a non-constant element $t\in A$ and we let 
$$\phi_t=\sum_{i=0}^{r}a_i\tau ^i.$$

The statement of Conjecture ~\ref{C:Con} is not affected if we replace $K$ by a 
finite extension $K'$ since if we find a constant $C'$ that works for $K'$ in 
Conjecture ~\ref{C:Con}, then $C=\frac{C'}{[K':K]}$ will work for $K$.

If we conjugate $\phi$ by $\gamma\in K^{\alg}\setminus\{0\}$ (i.e. $a\rightarrow\gamma^{-1}\phi_a\gamma$ for every $a\in A$), we obtain a new Drinfeld module, 
which we denote by $\phi^{(\gamma)}$ and these two Drinfeld modules are 
isomorphic over $K(\gamma)$. As a particular case of Proposition $2$ of 
\cite{Poo} we get that $\hhat_{\phi}(x)=\hhat_{\phi^{(\gamma)}}(\gamma^{-1}x)$. 
Then, if Conjecture ~\ref{C:Con} is proved for $\phi^{(\gamma)}$, it will also 
hold for $\phi$. So, having 
these in mind, we replace $\phi$ by one of its conjugates that has the property 
that $\phi^{(\gamma)}_t$ is monic, i.e. with the above 
notations, $\gamma$ satisfies the equation $\gamma^{q^r-1}a_r=1$. Because $[K(\gamma):K]\le q^r-1$, working over $K(\gamma)$ instead of $K$, we may introduce a factor of $\frac{1}{q^r-1}$ at the worst in the constant $C$ from Conjecture ~\ref{C:Con}, as explained in the previous paragraph. Also, the 
module structure theorems that we will be proving in the next section will 
not be affected by replacing $\phi$ with an isomorphic Drinfeld module or by replacing $K$ with a finite extension.

So, for simplifying the notation we suppose from now on in this section that $\phi_t$ is monic. 

In this section we will prove Theorem~\ref{T:T2mw'},
which 
is a special case of the Lehmer inequality for
Drinfeld modules. We will actually prove a more
general result (Theorem~\ref{T:T2mwg}) from which we
will be able to infer Theorems ~\ref{T:T2mw} and
\ref{T:T2mw'}.

For each
finite extension $L$ of $K$, we let $S_L$ be the set
of places $v\in M_L$ for which there exists a
coefficient $a_i$ of $\phi_t$ such that $v(a_i)<0$. We will prove that the set $S_L$ is the set of all valuations on $L$ of \emph{bad reduction} for $\phi$. We define next the notion of \emph{good reduction} for a Drinfeld module.

\begin{definition}
\label{D:badred}
Let $\phi:A\rightarrow K\{\tau\}$ be a Drinfeld
module. Let $L$ be a finite extension of $K$. We call
$v\in M_L$ a place of \emph{good} reduction for $\phi$
if for all $a\in A\setminus\{0\}$, the coefficients of
$\phi_a$ are integral at $v$ and the leading
coefficient of $\phi_a$ is a unit in the valuation
ring at $v$. If $v\in M_L$ is not a place of good
reduction, we call it a place of \emph{bad} reduction.
\end{definition}

\begin{lemma}
\label{L:bad reduction}
The set $S_L$ is the set of all places from $M_L$ at
which $\phi$ has bad reduction.
\end{lemma}

\begin{proof}
By the construction of the set $S_L$, the places from
$S_L$ are of bad reduction for $\phi$. We will prove
that these are all the bad places for $\phi$. 

Let $a\in A$. The equation $\phi_a\phi_t=\phi_t\phi_a$
will show 
that all the places where not all of the coefficients
of $\phi_a$ are integral, are from $S_L$. Suppose this
is 
not the case and take a place $v\notin S_L$ at which
some coefficient of $\phi_a$ is not integral. Let
$\phi_a=\sum_{i=0}^{r'}a_i'\tau^i$ and assume 
that $i$ is the largest index for a coefficient $a'_i$
that is not integral at
$v$. 

We equate the coefficient of $\tau^{i+r}$ in
$\phi_a\phi_t$ and $\phi_t\phi_a$, respectively. The
former is 
\begin{equation}
\label{E:e1}
a_i'+\sum_{j>i}a_j'a_{r+i-j}^{q^j}
\end{equation}
while the latter is 
\begin{equation}
\label{E:e2}
a_i'^{q^r}+\sum_{j>i}a_{r+i-j}a_j'^{q^{r+i-j}}.
\end{equation}
Thus the valuation at $v$ of \eqref{E:e1} is
$v(a_i')$, because all the $a_j'$ (for $j>i$) and
$a_{r+i-j}$ are integral at $v$, while $v(a_i')<0$.
Similarly, the valuation 
of \eqref{E:e2} is $v(a_i'^{q^r})=q^rv(a_i')<v(a_i')$
($r\ge 1$ because $t$ is non-constant). This fact
gives a 
contradiction to $\phi_a\phi_t=\phi_t\phi_a$. So, the
coefficients of $\phi_a$ for all $a\in A$, are
integral at all places of $M_L\setminus S_L$.

Now, using the same equation
$\phi_a\phi_t=\phi_t\phi_a$ and equating the leading
coefficients in both polynomials we obtain
$$a'_{r'}=a_{r'}'^{q^r}.$$
So, $a'_{r'}\in\F_p^{\alg}$. Thus, all the
leading coefficients for polynomials $\phi_a$ are
constants with respect to the valuations of $L$. So, if $v\in M_L\setminus S_L$, then all
the coefficients of $\phi_a$ are integral at $v$ and
the leading coefficient of $\phi_a$ is a unit in the
valuation ring at $v$ for every $a\in
A\setminus\{0\}$. Thus, $v\notin S_L$ is a place of
good reduction for $\phi$.
\end{proof}

\begin{theorem}
\label{T:T2mw}
Let $K$ be a finitely generated field of characteristic $p$. Let
$\phi:A\rightarrow K\{\tau\}$ be a Drinfeld module and
assume that there exists a non-constant $t\in A$ such
that $\phi_t$ is monic. Let $F$ be the algebraic
closure of $\F_p$ in $K$. We let $M_{K/F}$ be
the coherent good set of valuations on $K$, constructed as in
Section $3$. Let $\hhat$ and $\hhat_v$ be the global and local 
heights associated to $\phi$, constructed with respect to the coherent good set 
of valuations $M_{K/F}$. Let $x\in K^{\alg}$ and
let $F_x$ be the algebraic closure of $\F_p$
in $K(x)$. We construct the good set of valuations $M_{K(x)/F_x}$ which lie 
above the valuations from $M_{K/F}$. Let $S_x$ be the set of places $v\in 
M_{K(x)/F_x}$ such that $\phi$ has bad reduction at $v$. 

If $x$ is not a torsion point for $\phi$, then
there exists $v\in M_{K(x)/F_x}$ such that 
$$\hhat_v(x)>q^{-r\left(2+(r^2+r)|S_x|\right)}d(v)$$
where $d(v)$ is as always the degree of the valuation $v$.
\end{theorem}

Let $\{x_1,\dots,x_n\}$ be the transcendence basis for
$K/F$ associated to the construction of $M_{K/F}$. Let $v_0\in M_{K/F}$ be the 
place lying below the place $v$ from the conclusion of Theorem~\ref{T:T2mw}. 
Then $d(v)=\frac{d(v_0)f(v|v_0)}{[K(x):K]}$. Because $f(v|v_0)\ge 1$, $d(v_0)\ge
\frac{1}{[K:F(x_1,\dots,x_n)]}$ (see
\eqref{E:important inequality}) and $\hhat(x)\ge
\hhat_v(x)$, Theorem \ref{T:T2mw} has the following
corollary.
\begin{theorem}
\label{T:T2mw'}
With the notation from Theorem \ref{T:T2mw}, if
$x\notin\phi_{\tor}$, then
$$\hhat(x)>\frac{q^{-r\left(2+(r^2+r)|S_x|\right)}}{\left[
K(x):F(x_1, \dots,x_n)\right]}.$$
\end{theorem}

\begin{remark}
Theorem \ref{T:T2mw'} is a weaker form of Conjecture~\ref{C:Con} because our constant $C$ for which
$\hhat(x)\ge\frac{C}{[K(x):K]}$ for
$x\notin\phi_{\tor}$, is not completely
independent of $K(x)$. For us,
$$C=\frac{q^{-r\left(2+(r^2+r)|S_x|\right)}}{\left[K:F(x_1,\dots,
x_n\right)]}$$
and $S_x$ depends on $K(x)$.
\end{remark}

Before proving Theorem~\ref{T:T2mw}, we need to prove several preliminary lemmas regarding the local height for an arbitrary point of the Drinfeld module $\phi$. 

Fix now a finite extension $L$ of $K$ and let $U$ be a
good set of valuations on $L$. Let $S=S_L\cap U$.

For each $v\in U$ we define
\begin{equation}
\label{E:defM_v2}
M_v=\min_{i\in\{0,\dots,r-1\}}\frac{v(a_i)}{q^r-q^i}
\end{equation}
where by convention, as always: $v(0)=+\infty$. We observe that
$M_v<0$ if and only if $v\in S$.

Let $v\in S$. We define $P_v$ as the subset of the negatives of the slopes
of the Newton polygon associated to $\phi_t$, consisting of those $\alpha$ for 
which there exist $i\ne j$ in $\{0,\dots,r\}$ such that
\begin{equation}
\label{E:defP_v2}
\alpha=\frac{v(a_i)-v(a_j)}{q^j-q^i}\le 0,
\end{equation}
and $v(a_i)+q^i\alpha=v(a_j)+q^j\alpha=\min_{0\le l\le 
r}\left(v(a_l)+q^l\alpha\right)$. If $\phi$ is the Carlitz module in 
characteristic $2$,
i.e. $\phi=\psi_2$, where
$\psi_2:\F_2[t]\rightarrow K\{\tau\}$ is
defined by $\psi_2(x)=tx+x^2$ for every $x$, then we
want the set $P_v$ to contain $\{0\}$, even if $0$ is
not in the set from \eqref{E:defP_v2}.

Let \[ N_{\phi}=\begin{cases}
1+r=2, \text{ if }
\phi=\psi_2\\
r, \text{ otherwise.}
\end{cases}
\] 

Clearly, for every $\phi$ and $v\in S$, $|P_v|\le
N_{\phi}$. We define next the concept of angular component for a nonzero $x\in L$. For this we first fix a uniformizer $\pi_v\in L$ for each valuation $v\in S$.

\begin{definition}
\label{D:angular}
Assume $v\in S$. For every nonzero $y\in L$ we define
the angular component of 
$y$ at $v$, denoted by $\ac_{\pi_v}(y)$, to be the residue
at $v$ of $y\pi_v^{-v(y)}$. 
(Note that the angular component is never $0$.)
\end{definition}
We can define in a similar manner as above the notion
of angular component at 
each $v\in M_L$ but we will work with angular
components at the places from $S$ 
only.

The main property of the angular component is that for
every $y,z\in L\setminus\{0\}$, 
$v(y-z)>v(y)=v(z)$ if and only if
$(v(y),\ac_{\pi_v}(y))=(v(z),\ac_{\pi_v}(z))$.

For each $\alpha\in P_v$ we let $l\ge 
1$ and let $ i_0<i_1<\dots <i_l$ be all the indices
$i$ for which 
\begin{equation}
\label{E:Newton2}
v(a_i\alpha^{q^i})=\min_jv(a_j\alpha^{q^j}).
\end{equation}
We define $R_v(\alpha)$ as the set containing all the nonzero solutions of the equation 

\begin{equation}
\label{E:defR_v2}
\sum_{j=0}^{l}\ac_{\pi_v}(a_{i_j})X^{q^{i_j}}=0,
\end{equation}
where the indices $i_j$ are the ones associated to $\alpha$ as in \eqref{E:Newton2}.
For $\alpha=0$, we want the set $R_v(\alpha)$ to
contain also $\{1\}$ in addition to the numbers from
\eqref{E:defR_v2}. If $\alpha=0$, $l$ might be $0$ and so, equation \eqref{E:defR_v2} might have no nonzero solutions. In that case, $R_v(0)=\{1\}$. Clearly, for every $\alpha\in
P_v$, $|R_v(\alpha)|\le q^r$. 

\begin{lemma}
\label{L:L0}
Assume $v\in S$ and let $x\in L$. If $v(x)\le 0$ and
$v(\phi_t(x))>\min_{i\in\{0,\dots,r\}}v(a_ix^{q^i})$
then 
$(v(x),\ac_{\pi_v}(x))\in P_v\times R_v(v(x))$.
\end{lemma}
\begin{proof}
If $v(\phi_t(x))>\min_{i\in\{0,\dots,r\}}v(a_ix^{q^i})$ it
means that there exists 
$l\ge 1$ and $$i_0<\dots <i_l$$ such that
\begin{equation}
\label{E:0,l,v}
v(a_{i_0}x^{q^{i_0}})=\dots 
=v(a_{i_l}x^{q^{i_l}})=\min_{i\in\{0,\dots,r\}}v(a_ix^{q^i})
\end{equation} 
and also 
\begin{equation}
\label{E:0,l,ac}
\sum_{j=0}^{l}\ac_{\pi_v}(a_{i_j})\ac_{\pi_v}(x)^{q^{i_j}}=0.
\end{equation}
Equations \eqref{E:0,l,v} and \eqref{E:0,l,ac} yield
$v(x)\in P_v$ and 
$\ac_{\pi_v}(x)\in R_v(v(x))$ respectively, according
to \eqref{E:defP_v2} and \eqref{E:defR_v2}. 

Note that we needed our assumption that $v(x)\le 0$ only because $P_v$ consists only of the negatives of the non-negative slopes of the Newton polygon associated to $\phi_t$ (and not the negatives of all the slopes). The above proof shows that as long as the valuation of $x$ and the angular component of $x$ do not belong to certain prescribed sets, $v(\phi_t(x))=\min_iv(a_ix^{q^i})$.
\end{proof}

\begin{lemma}
\label{L:L2'}
Let $v\in M_L$ and let $x\in L$. If
$v(x)<\min\{0,M_v\}$, then 
$\hhat_v(x)=-d(v)\cdot v(x)$.
\end{lemma}

\begin{proof}  
For every $i\in\{0,\dots,r-1\}$,
$v(a_ix^{q^i})=v(a_i)+q^iv(x)>q^rv(x)$ 
because
$v(x)<M_v=\min_{i\in\{0,\dots,r-1\}}\frac{v(a_i)}{q^r-q^i}$.
This 
shows that $v(\phi_t(x))=q^rv(x)<v(x)<\min\{0,M_v\}$.
By induction, 
$v(\phi_{t^n}(x))=q^{rn}v(x)$ 
for all $n\ge 1$. So, $V_v(x)=v(x)$ 
and 
$$\hhat_v(x)=-d(v)\cdot v(x).$$
\end{proof}

An immediate corollary to \eqref{L:L2'} is the
following result.
\begin{lemma}
\label{L:L1'}
Assume $v\notin S$ and let $x\in L$. If $v(x)<0$ then
$\hhat_v(x)=-d(v)\cdot v(x)$, 
while if $v(x)\ge 0$ then $\hhat_v(x)=0$.
\end{lemma}
\begin{proof}  
First, it is clear that if $v(x)\ge 0$ then for all $
n\ge 1$, 
$v(\phi_{t^n}(x))\ge 0$ because all the coefficients
of $\phi_t$ and thus of 
$\phi_{t^n}$ have non-negative valuation at $v$. Thus
$V_v(x)=0$ and so, 
$$\hhat_v(x)=0.$$
Now, if $v(x)<0$, then $v(x)<M_v$ because $M_v\ge 0$
($v\notin S$). So, applying 
the result of \eqref{L:L2'} we conclude the proof of
this lemma.
\end{proof}

We will get a better insight into the local heights
behaviour with the following 
lemma.
\begin{lemma}
\label{L:L3'}
Let $x\in L$. Assume $v\in S$ and $v(x)\le 0$. If
$(v(x),\ac_{\pi_v}(x))\notin P_v\times 
R_v(v(x))$ then $v(\phi_t(x))<M_v$, unless $q=2$,
$r=1$ and $v(x)=0$.
\end{lemma}
\begin{proof}
Lemma \ref{L:L0} implies that there exists
$i_0\in\{0,\dots,r\}$ such that 
for all $i\in\{0,\dots,r\}$ we have 
$v(a_ix^{q^i})\ge v(a_{i_0}x^{q^{i_0}})=v(\phi_t(x))$.

Suppose \eqref{L:L3'} is not true and so, there exists
$j_0<r$ such that
$$\frac{v(a_{j_0})}{q^r-q^{j_0}}\le
v(\phi_t(x))=v(a_{i_0})+q^{i_0}v(x).$$
This means that
\begin{equation}
\label{E:a_j_0_1}
v(a_{j_0})\le
(q^r-q^{j_0})v(a_{i_0})+(q^{r+i_0}-q^{i_0+j_0})v(x).
\end{equation}
On the other hand, by our assumption about $i_0$, we
know that 
$v(a_{j_0}x^{q^{j_0}})\ge v(a_{i_0}x^{q^{i_0}})$ which
means that 
\begin{equation}
\label{E:a_j_0_2}
v(a_{j_0})\ge v(a_{i_0})+(q^{i_0}-q^{j_0})v(x).
\end{equation}
Putting together inequalities \eqref{E:a_j_0_1} and
\eqref{E:a_j_0_2}, we get 
$$v(a_{i_0})+(q^{i_0}-q^{j_0})v(x)\le 
(q^r-q^{j_0})v(a_{i_0})+(q^{r+i_0}-q^{i_0+j_0})v(x).$$
Thus 
\begin{equation}
\label{E:1'}
v(x)(q^{r+i_0}-q^{i_0+j_0}-q^{i_0}+q^{j_0})\ge
-v(a_{i_0})(q^r-q^{j_0}-1).
\end{equation}
But
$q^{r+i_0}-q^{i_0+j_0}-q^{i_0}+q^{j_0}=q^{r+i_0}(1-q^{j_0-r}-q^{-r}+q^{j_0-r-i_0
})$ and because $j_0<r$ and $q^{j_0-r-i_0}>0$, we
obtain
\begin{equation}
\label{E:2'}
1-q^{j_0-r}-q^{-r}+q^{j_0-r-i_0}>1-q^{-1}-q^{-r}\ge
1-2q^{-1}\ge 0.
\end{equation}
Also, $q^r-q^{j_0}-1\ge
q^r-q^{r-1}-1=q^{r-1}(q-1)-1\ge 0$ with equality if
and 
only if 
$q=2$, $r=1$ and $j_0=0$. We will analyze this case
separately. So, as long as 
we are not in this special case, we do have 
\begin{equation}
\label{E:3'}
q^r-q^{j_0}-1>0.
\end{equation}
Now we have two possibilities (remember that $v(x)\le
0$):

(i)  $v(x)<0$  

In this case, \eqref{E:1'}, \eqref{E:2'} and
\eqref{E:3'} tell us that 
$-v(a_{i_0})<0$. Thus, $v(a_{i_0})>0$. But we know
from our hypothesis on 
$i_0$ that $v(a_{i_0}x^{q^{i_0}})\le v(x^{q^r})$ which
is in contradiction with 
the combination of the following facts: $v(x)<0$,
$i_0\le r$ and $v(a_{i_0})>0$.

(ii)  $v(x)=0$

Then another use of \eqref{E:1'}, \eqref{E:2'} and
\eqref{E:3'} gives us 
$-v(a_{i_0})\le 0$; thus $v(a_{i_0})\ge 0$. This would
mean that 
$v(a_{i_0}x^{q^{i_0}})\ge 0$ and this contradicts our
choice for $i_0$ because 
we know from the fact that $v\in S$, that there exists
$i\in \{0,\dots,r\}$ 
such that $v(a_i)<0$. So, then we would have 
$$v(a_ix^{q^i})=v(a_i)<0\le v(a_{i_0}x^{q^{i_0}}).$$
Thus, in either case (i) or (ii) we get a
contradiction that proves the lemma 
except in the 
special case that we excluded above: $q=2$, $r=1$ and
$j_0=0$. If we have $q=2$ 
and $r=1$ then 
$$\phi_t(x)=a_0x+x^2.$$
By the definition of $S$ and because $v\in S$, $v(a_0)<0$. Also, $M_v=v(a_0)$.

If $v(x)<0$, then either $v(x)<M_v=v(a_0)$, in
which case again 
$v(\phi_t(x))<M_v$ (as shown in the proof of Lemma~\ref{L:L2'}), or $v(x)\ge 
M_v$. In the latter case, 
$$v(\phi_t(x))=v(a_0x)=v(a_0)+v(x)<v(a_0)=M_v.$$
So, we see that indeed, only
$v(x)=0$, $q=2$ and $r=1$ can make 
$v(\phi_t(x))\ge M_v$ in the hypothesis of
\eqref{L:L3'}. 
\end{proof}

\begin{lemma}
\label{L:L4'}
Assume $v\in S$ and let $x\in L$. Excluding the case $q=2$, $r=1$ and
$v(x)=0$, we have that if 
$v(x)\le 0$ 
then either $\hhat_v(x)>\frac{-d(v)M_v}{q^r}$ or 
$(v(x),ac_{\pi_v}(x))\in P_v\times R_v(v(x))$.
\end{lemma}

\begin{proof}
If $v(x)\le 0$ then

$$\text{\emph{either}: (i) } v(\phi_t(x))<M_v\text{
,}$$

in which case by \eqref{L:L2'} we have that 
$\hhat_v(\phi_t(x))=-d(v)\cdot v(\phi_t(x))$. So, case (i) yields
\begin{equation}
\label{E:4'}
\hhat_v(x)=-d(v)\cdot\frac{v(\phi_t(x))}{\deg\phi_t} > 
-d(v)\cdot\frac{M_v}{q^r}
\end{equation}

$$\text{\emph{or}: (ii) }v(\phi_t(x))\ge M_v\text{
,}$$

in which case, Lemma~\ref{L:L3'} yields 
\begin{equation}
\label{E:(ii2)}
v(\phi_t(x))>v(a_{i_0}x^{q^{i_0}})=\min_{i\in\{0\dots,r\}}v(a_ix^{q^i}).
\end{equation}
Using \eqref{E:(ii2)} and Lemma~\ref{L:L0} we
conclude that case (ii) yields 
$(v(x),\ac_{\pi_v}(x))\in P_v\times R_v(v(x))$.
\end{proof}

Now we analyze the excluded case from Lemma~\ref{L:L4'}.
\begin{lemma}
\label{L:L5'}
Assume $v\in S$ and let $x\in L$. If $v(x)\le 0$, then either
$$(v(x),ac_{\pi_v}(x))\in P_v\times 
R_v(v(x))$$ or $\hhat_v(x)\ge\frac{-d(v)M_v}{q^r}$.
\end{lemma}

\begin{proof}
Using the result of \eqref{L:L4'} we have left to
analyze the case: $q=2$, $r=1$ 
and $v(x)=0$. 

As shown in the proof of \eqref{L:L3'}, in this case
$\phi_t(x)=a_0x+x^2$ and 
$$v(\phi_t(x))=v(a_0)=M_v<0.$$ 
Then, either
$v(\phi_{t^2}(x))=v(\phi_t(x)^2)=2M_v<M_v$ or 
$v(\phi_{t^2}(x))>v(a_0\phi_t(x))=v(\phi_t(x)^2)$. If
the former case holds, 
then by \eqref{L:L2'},
$$\hhat_v(\phi_{t^2}(x))=-d(v)\cdot 
2M_v$$ and so,
$$\hhat_v(x)=\frac{-d(v)\cdot 2M_v}{4}.$$
If the latter case holds, i.e. 
$v(\phi_{t}(\phi_t(x)))>v(a_0\phi_t(x))=v(\phi_t(x)^2)$,
then 
$\ac_{\pi_v}(\phi_t(x))$ satisfies the equation
$$\ac_{\pi_v}(a_0)X+X^2=0.$$ 
Because the angular component is never $0$, it must be
that 
$\ac_{\pi_v}(\phi_t(x))=\ac_{\pi_v}(a_0)$ (remember
that we are working now in 
characteristic $2$). But, because 
$v(a_0x)<v(x^2)$ we can relate the angular component
of $x$ and the angular 
component of $\phi_t(x)$ and so,
$$\ac_{\pi_v}(a_0)=\ac_{\pi_v}(\phi_t(x))=\ac_{\pi_v}(a_0x)=
\ac_{\pi_v}(a_0)\ac_{\pi_v}(x).$$
This means $\ac_{\pi_v}(x)=1$ and so, the excluded
case amounts to a 
dichotomy similar to the one from \eqref{L:L4'}:
either 
$(v(x),\ac_{\pi_v}(x))=(0,1)$ or
$\hhat_v(x)=\frac{-d(v)M_v}{2}$. The definitions of $P_v$
and $R_v(\alpha)$ from 
\eqref{E:defP_v2} and \eqref{E:defR_v2} respectively,
yield that $(0,1)\in 
P_v\times R_v(0)$.
\end{proof}

Finally, we note that in \eqref{L:L5'} we have 
$$-\frac{d(v)M_v}{q^r}=-\frac{d(v)e(v|v_0)M_{v_0}}{q^r}.$$
We have obtained the following dichotomy.
 
\begin{lemma}
\label{L:L5}
Assume $v\in S$. If $v(x)\le 0$ then either
$(v(x),ac_{\pi_v}(x))\in 
P_v\times R_v(v(x))$ or
$\hhat_{U,v}(x)\ge\frac{-M_vd(v)}{q^r}$. Moreover, by
its definition $M_v<-\frac{1}{q^r}$ and so, if the above latter case holds for 
$x$, then
$\hhat_{U,v}(x)>\frac{d(v)}{q^{2r}}$.
\end{lemma}

We will deduce Theorem \ref{T:T2mw} from the following
more general result.
\begin{theorem}
\label{T:T2mwg}
Let $K$ be a field extension of $\F_q$ and let
$\phi:A\rightarrow K\{\tau\}$ be a Drinfeld module.
Let $L$ be a finite field extension of $K$. Let $t$ be
a non-constant element of $A$ and assume that
$\phi_t=\sum_{i=0}^ra_i\tau^i$ is monic. Let $U$ be a
good set of valuations on $L$ and let $C(U)$ be, as always, the field of 
constants with respect to $U$. Let $S$ be the finite
set of valuations $v\in U$ such that $\phi$ has bad reduction at $v$.
Let $x\in L$.

a) If $S$ is empty, then either $x\in C(U)$ or there
exists $v\in U$ such that $\hhat_{U,v}(x)\ge d(v)$.

b) If $S$ is not empty, then either $x\in\phi_{\tor}$,
or there exists $v\in U$ such that
$\hhat_{U,v}(x)>q^{-2r-r^2N_{\phi}|S|}d(v)\ge
q^{-r\left(2+(r^2+r)|S|\right)}d(v)$.
Moreover, if $S$ is not empty and $x\in\phi_{\tor}$,
then there exists a polynomial
$b(t)\in\F_q[t]$ of degree at most
$rN_{\phi}|S|$ such that $\phi_{b(t)}(x)=0$.
\end{theorem}

\begin{proof} 
$a)$ Assume $S$ is empty.

Then either $v(x)\ge 0$ for all $v\in U$ or there
exists $v\in U$ such that $v(x)<0$. If the latter
statement is true, then Lemma~\ref{L:L1'} shows that
for any valuation $v\in U$ for which $v(x)<0$, we have
$$\hhat_{U,v}\ge d(v)\text{,}$$
because $v\notin S$ ($S$ is empty).

Now, if the former statement is true, i.e. $x$ is
integral at all places from $U$, then by Lemma
~\ref{L:Lconstants}, $x\in C(U)$.

$b)$ Assume $S$ is not empty.

Let $v\in S$. We will use several times the following result.
\begin{lemma}
\label{L:L-1}
Let $L$ be a field extension of $\F_q$ and let $v$ be a discrete 
valuation on 
$L$. Let $I$ be a finite set of integers. Let $N$ be an
integer greater or equal than 
all the elements of $I$. For each $\alpha\in I$, let
$R(\alpha)$ be a nonempty 
finite set 
of nonzero elements of the residue field at $v$. Let
$W$ be an 
$\F_q$-vector subspace of $L$ with the
property that for all $w\in 
W$, $(v(w),\ac_{\pi_v}(w))\in I\times R(v(w))$
whenever $v(w)\le N$.

Let $f$ be the smallest positive integer greater or equal than
$\max_{\alpha\in 
I}\log_{q}\vert R(\alpha)\vert$. Then the $\F_q$-codimension
of $\left\{w\in W\mid 
v(w)>N\right\}$ is bounded above by $\vert I\vert f$.
\end{lemma}

\begin{proof}[Proof of Lemma~\ref{L:L-1}.]
Let $i=|I|$. Let $\alpha_0<\dots <\alpha_{i-1}$ be the elements of $I$, and let 
$\alpha_i=N+1$. For $0\le j\le i$, define $W_j=\{w\in W|v(w)\ge\alpha_j\}$. For 
$0\le j<i$, the hypothesis gives an injection 
$$W_j/W_{j+1}\rightarrow R(\alpha_j)\cup\{0\}$$
taking $w$ to the residue of $w/\pi_v^{\alpha_j}$. Thus 
$$q^{\dim_{\F_q}W_j/W_{j+1}}\le q^f+1<q^{f+1},$$
so $\dim_{\F_q}W_j/W_{j+1}\le f$ (note that we used the fact that $f>0$ 
in order to have the inequality $q^f+1<q^{f+1}$). Summing over $j$ gives 
$\dim_{\F_q}W_0/W_i\le if$, as desired.
\end{proof}

We apply Lemma~\ref{L:L-1} with $N=0$,
$I=P_v$ and $R(\alpha)=R_v(\alpha)$ for every
$\alpha\in P_v$. Because $|P_v|\le N_{\phi}$ and
$|R_v(\alpha)|\le q^r$ for every $\alpha\in P_v$, we
obtain the following result.

\begin{Fact}
\label{C:Cor0}
Let $v\in S$. Let $W$ be an $\F_q$-subspace of
$L$ with the property that for all $w\in W$,
$(v(w),\ac_{\pi_v}(w))\in P_v\times R_v(v(w))$
whenever $v(w)\le 0$. 

Then the $\F_q$-codimension of $\left\{w\in W\mid
v(w)>0\right\}$ in $W$ is bounded above by $rN_{\phi}$.
\end{Fact}

We apply Fact~\ref{C:Cor0} for each $v\in S$ and
we deduce the following two results.
\begin{Fact}
\label{C:Cor1}
Let $W$ be an $\F_q$-subspace of $L$ such that
$(v(w),\ac_{\pi_v}(w))\in P_v\times R_v(v(x))$
whenever $v\in S$, $w\in W$ and $v(w)\le 0$. Then the
$\F_q$-codimension of 
$$\left\{w\in W\mid v(w)>0\text{ for all $v\in
S$}\right\}$$
in $W$ is bounded above by $rN_{\phi}|S|$.
\end{Fact}

\begin{Fact}
\label{C:Cor2}
Let notation be as in Fact~\ref{C:Cor1}. If
moreover, $\dim_{\F_q}W>rN_{\phi}|S|$, then
there exists a nonzero $w\in W$ such that $v(w)>0$ for
all $v\in S$.
\end{Fact}

Using the above facts we prove the following claim which is the key for 
obtaining the result of Theorem~\ref{T:T2mwg}.
\begin{Claim}
\label{L:L-2}
Assume $|S|\ge 1$. If $W$ is an
$\F_q$-subspace of $L$ and
$\dim_{\F_q}W>rN_{\phi}|S|$, then there exists
$w\in W$ and there exists $v\in U$ such that
$\hhat_{U,v}(w)>\frac{d(v)}{q^{2r}}$.
\end{Claim}

\begin{proof}[Proof of Claim \ref{L:L-2}.]
If there exists $v\in U\setminus S$ and $w\in W$ such
that $v(w)<0$, then by Lemma~\ref{L:L1'},
$$\hhat_v(w)\ge d(v)>\frac{d(v)}{q^{2r}}.$$

Thus, suppose from now on in the proof of Claim~\ref{L:L-2},
that for every $v\in U\setminus S$ and every $w\in W$,
$v(w)\ge 0$.

Because $\dim_{\F_q}W>rN_{\phi}|S|$, Fact~\ref{C:Cor2} guarantees the 
existence of a nonzero
$w\in W$ for which \emph{either} $v(w)>0$ for all
$v\in S$, \emph{or} there exists $v\in S$ such that 
\begin{equation}
\label{E:contr}
v(w)\le 0\text{ but }(v(w),\ac_{\pi_v}(w))\notin
P_v\times R_v(v(w)).
\end{equation}
The former case is impossible because we already
supposed that $v(w)\ge 0$ for all $v\in U\setminus S$.
Because $|S|\ge 1$ there is no nonzero $w$ that has
non-negative valuation at all the places in $U$ and
positive valuation at at least one place in $U$. Its
existence would contradict the sum formula for the
valuations in $U$.

Thus, the latter case holds, i.e. there exists $v\in
S$ satisfying \eqref{E:contr}. But then, Lemma
~\ref{L:L5} gives
$\hhat_{U,v}(w)>\frac{d(v)}{q^{2r}}$.
\end{proof}

Using Claim \ref{L:L-2} we can finish the proof of
part $b)$ of Theorem~\ref{T:T2mwg}.

Consider $W=\Span_{\F_q}\left(\left\{
x,\phi_t(x),\dots,\phi_{t^{rN_{\phi}|S|}}(x)\right\}\right)$.
If there exists no polynomial $b(t)$ as in the
statement of part $b)$ of Theorem~\ref{T:T2mwg}, then
$\dim_{\F_q}W=1+rN_{\phi}|S|$. Applying Claim~\ref{L:L-2} to $W$, we 
find some $w\in W$ and some
$v\in U$ such that 
\begin{equation}
\label{E:wineq}
\hhat_{U,v}(w)>\frac{d(v)}{q^{2r}}.
\end{equation}
By the construction of $W$, then there exists
a nonzero polynomial $d(t)\in\F_q[t]$ of
degree at most $rN_{\phi}|S|$ such that 
\begin{equation}
\label{E:wx}
w=\phi_{d(t)}(x).
\end{equation}
Using equations \eqref{E:wineq} and \eqref{E:wx}, we
obtain
$$\hhat_{U,v}(x)=\frac{\hhat_{U,v}(w)}{\deg(\phi_{d(t)})}>\frac{\frac{d(v)}{q^{2
r}}}
{q^{r\cdot rN_{\phi}|S|}}\text{ ,}$$
as desired.
\end{proof}

\begin{proof}[Proof of Theorem \ref{T:T2mw}.]
There are two cases.

\emph{Case 1.} The set $S_x$ is empty.

By Lemma \ref{L:Lconstants}, all the coefficients
$a_i$ of $\phi_t$ are from $F_x$. Let
$\F_{q^l}$ be a finite field containing all
these coefficients.

Assume $x\in\F_p^{\alg}$. Let
$\F_{q^{ll'}}=\F_{q^l}(x)$. Then for
every $n\ge 1$,
$\phi_{t^n}(x)\in\F_{q^{ll'}}$. Because
$\F_{q^{ll'}}$ is finite, there exist distinct
positive integers $n$ and $n'$ such that
$\phi_{t^n}(x)=\phi_{t^{n'}}(x)$. Thus
$\phi_{t^{n'}-t^n}(x)=0$; i.e. $x\in\phi_{\tor}$,
which is a contradiction with our hypothesis that $x$
is not a torsion point.

Thus, in \emph{Case 1}, $x\notin\F_p^{\alg}$.
So, $x$ is not a constant with respect to the
valuations from $M_{K(x)/F_x}$. Then, by Theorem
~\ref{T:T2mwg} $a)$, there exists $v\in M_{K(x)/F_x}$
such that 
$$\hhat_v(x)\ge d(v)>q^{-2r}d(v).$$

\emph{Case 2.} The set $S_x$ is not empty.

Because $x\notin\phi_{\tor}$, Theorem \ref{T:T2mwg}
shows that there exists $v\in M_{K(x)/F_x}$ such that 
$$\hhat_v(x)> q^{-2r-r^2N_{\phi}|S_x|}d(v)\ge
q^{-r\left(2+(r^2+r)|S_x|\right)}d(v).$$
\end{proof}

\begin{remark}
Assume that we have a Drinfeld module
$\phi:A\rightarrow K\{\tau\}$ and a
non-constant element $t\in A$ for which $\phi_t$ is monic. Suppose we are in 
\emph{Case 1} of the proof of Theorem~\ref{T:T2mw}. Then that proof
shows that for every non-torsion $x\in K^{\alg}$,
there exists $v\in M_{K(x)/F_x}$ such that
$\hhat_v(x)\ge\frac{d(v_0)}{[K(x):K]}$, where $v_0$ is
the place of $M_{K/F}$ that sits below $v$. Because of
inequality \eqref{E:important inequality},
$d(v_0)\ge\frac{1}{\left[K:F(x_1,\dots,x_n)\right]}$,
where $\{x_1,\dots,x_n\}$ is the transcendence basis
for $K/F$ with respect to which we constructed the
good set of valuations $M_{K/F}$. Thus Conjecture~\ref{C:Con} holds in this case, i.e. when all the
coefficients $a_i$ are from $\F_p^{\alg}$,
with $C=\frac{1}{\left[K:F(x_1,\dots,x_n)\right]}$.
\end{remark}

With the help of Theorem \ref{T:T2mw} we can get a
characterization of the torsion submodule of a
Drinfeld module. Let $K$ be a finitely generated field
and let $\phi:A\rightarrow K\{\tau\}$ be a Drinfeld
module. If none of the non-constant $a\in A$ has the
property that $\phi_a$ is monic, then just pick some
non-constant $t\in A$ and conjugate $\phi$ by
$\gamma\in K^{\alg}\setminus\{0\}$ such that $\phi^{(\gamma)}_t$ is
monic. Then $\phi$ and $\phi^{(\gamma)}$ are
isomorphic over $K(\gamma)$, which is a finite
extension of $K$ of degree at most $\deg(\phi_t)-1$.
So, describing $\phi_{\tor}\left(K(\gamma)\right)$ is
equivalent with describing
$\phi^{(\gamma)}_{\tor}\left(K(\gamma)\right)$. The
following result does exactly this. Its proof is
immediate after the proof of Theorem~\ref{T:T2mwg}.

\begin{corollary}
\label{C:CorT2}
Let $K$ be a finitely generated field and let
$\phi:A\rightarrow K\{\tau\}$ be a Drinfeld module.
Let $t$ be a non-constant element of $A$. Let
$\phi_t=\sum_{i=0}^{r}a_i\tau^i$ and assume that
$a_r=1$. Let $L$ be a finite extension of $K$ and let
$E$ be the algebraic closure of $\F_p$ in $L$.

a) If $a_0,\dots,a_{r-1}\in E$, then
$\phi_{\tor}(L)=E$.

b) If not all of the coefficients $a_0,\dots,a_{r-1}$
are in $E$, let $S=S_L\cap M_{L/E}$. Let
$b(t)\in\F_q[t]$ be the least common multiple
of all the polynomials of degree at most
$rN_{\phi}|S|$. Then for all $x\in\phi_{\tor}(L)$,
$\phi_{b(t)}(x)=0$.
\end{corollary}

\begin{remark}
\label{R:conexiune}
We can also bound the size of the torsion of a Drinfeld module $\phi$ over a fixed field $K$ by specializing $\phi$ at a place of good reduction. This is the classical method used to bound torsion for abelian varieties. The bound that we would obtain by using this more classical method will be much larger than the one from Corollary~\ref{C:CorT2} if $K$ contains a large finite field. However, because our bound is obtained through completely different methods, one can use both methods and then choose the better bound provided by either one.
\end{remark}

The bound from Corollary \ref{C:CorT2} $b)$ for the
torsion subgroup of $\phi(L)$ is sharp when $\phi$ is
the Carlitz module, as shown by the following example.

\begin{example}
For each prime number $p$ let 
$v_{\infty}:\F_p(t)\setminus\{0\}\rightarrow\mathbb{Z}$ be the valuation such that 
$v(b)=-\deg(b)$ for each $b\in\F_p[t]\setminus\{0\}$. It is the same 
notation that
we used in Section $2$. Also, for each prime number
$p$, let $\psi_p$ be the Carlitz module in
characteristic $p$, i.e.
$\psi_p:\F_p[t]\rightarrow
\F_p(t)\{\tau\}$, given by
$\left(\psi_p\right)_t=t\tau^0+\tau$.

If $p=2$, we let $L=\F_2(t)$. Then with the
notation from Corollary ~\ref{C:CorT2},
$S=\{v_{\infty}\}$. Also, $r=1$,
$N_{\psi_2}=2$ and so,
$rN_{\psi_2}|S|=2$. It is immediate to see that
$\psi_2[t]\subset L$ and also $\psi_2[1+t]\subset L$.
Thus we do need a polynomial $b(t)$ of degree $2$, i.e.
$b(t)=t^2+t$, to kill the torsion of $\psi_2(L)$.

If $p>2$, we let
$L=\F_2\left((-t)^{\frac{1}{p-1}}\right)$.
Then $\psi_p[t]\subset L$. With the notation from
Corollary ~\ref{C:CorT2}, $r=1$ and $N_{\psi_p}=1$.
Also, $S=\{w_{\infty}\}$, where $w_{\infty}$ is the
unique place of $L$ sitting above $v_{\infty}$. So,
again we see that we need a polynomial $b(t)$ of
degree $rN_{\psi_p}|S|=1$ to kill the torsion of
$\psi_p(L)$.
\end{example}

\section{The Mordell-Weil theorem for infinitely generated fields}
Before stating and proving the theorems from this section we will introduce the 
notion of modular transcendence degree. This notion refers to the minimal field 
of definition for a Drinfeld module.

\begin{definition}
\label{D:fieldofdef}
For a Drinfeld module $\phi:A\rightarrow K\{\tau\}$, its field of definition is the smallest subfield of $K$ containing all the coefficients of $\phi_a$, for every $a\in A$.
\end{definition}

\begin{lemma}
\label{L:fieldofdeffingen}
The field of definition of a Drinfeld module is finitely generated.
\end{lemma}

\begin{proof}
Let $\phi:A\rightarrow K\{\tau\}$. Let $t_1,\dots,t_s$ be generators of $A$ as an $\F_q$-algebra. Let $K_0$ be the field extension of $\F_q$ generated by all the coefficients of $\phi_{t_1},\dots,\phi_{t_s}$. Then $K_0$ is finitely generated and by construction, $K_0$ is the field of definition for $\phi$.
\end{proof}

\begin{definition}
\label{D:modtrdeg}
Let $\phi:A\rightarrow K\{\tau\}$ be a Drinfeld
module. The modular 
transcendence degree of $\phi$ is the minimum
transcendence degree over 
$\F_p$ of the field of definition for
$\phi^{(\gamma)}$, where the 
minimum is taken over all $\gamma\in
K^{\alg}\setminus\{0\}$. 
\end{definition}

\begin{lemma}
\label{L:fieldef}
Let $\phi:A\rightarrow K\{\tau\}$ be a Drinfeld module
and let $E$ be its field of definition. Let $t\in A$
be a non-constant element and let
$\phi_t=\sum_{i=0}^{r}a_i\tau^i$. Let
$E_0=\F_p(a_0,\dots,a_r)$ and let $E_0^{\alg}$
be the algebraic closure of $E_0$ inside $K^{\alg}$.
Then $E_0\subset E\subset E_0^{\alg}$.
\end{lemma}

\begin{proof}
Let $\psi$ be the restriction of $\phi$ to $\F_p[t]$. Clearly, $\psi$ is defined over $E_0$. For every $a\in A$, $\phi_a$ is an endomorphism of $\psi$. Thus for every $a\in A$, by Proposition $4.7.4$ of \cite{Goss}, the coefficients of $\phi_a$ are algebraic over $E_0$.
\end{proof}

\begin{lemma}
\label{L:realmodtrdeg}
Let $\phi:A\rightarrow K\{\tau\}$ be a Drinfeld
module. Assume that there exists a non-constant
element $t\in A$ for which $\phi_t$ is monic. Let $E$
be the field of definition for $\phi$. Then the
modular transcendence degree of $\phi$ is
$\trdeg_{\F_p}E$.
\end{lemma}
\begin{proof}
By the definition of modular transcendence degree of
$\phi$, we have to show that for every $\gamma\in
K^{\alg}\setminus\{0\}$, if $E^{(\gamma)}$ is the field of
definition for $\phi^{(\gamma)}$, then 
\begin{equation}
\label{E:trgamma}
\trdeg_{\F_p}E^{(\gamma)}\ge\trdeg_{\F_p}E.
\end{equation}

Let $\gamma\in K^{\alg}\setminus\{0\}$. If
$\phi_t=\sum_{i=0}^r a_i\tau^i$, then
$\phi^{(\gamma)}_t=\sum_{i=0}^r
a_i\gamma^{q^i-1}\tau^i$. 

By Lemma~\ref{L:fieldef},
$\trdeg_{\F_p}E=\trdeg_{\F_p}\F_p(a_0,\dots,a_{r-1})$
and 
$$\trdeg_{\F_p}E^{(\gamma)}=\trdeg_{\F_p}\F_p\left
(a_0,a_1\gamma^{q-1},\dots,a_{r-1}\gamma^{q^{r-1}-1},\gamma^{q^r-1}\right).$$
So, in order to prove inequality \eqref{E:trgamma}, it
suffices to show that 
\begin{equation}
\label{E:t1}
\trdeg_{\F_p}\F_p\left
(a_0,a_1\gamma^{q-1},\dots,a_{r-1}\gamma^{q^{r-1}-1},\gamma^{q^r-1}\right)\ge
\trdeg_{\F_p}\F_p(a_0,\dots,a_{r-1}).
\end{equation}
But,
\begin{equation}
\label{E:t2}
\trdeg_{\F_p}\F_p\left
(a_0,\dots,a_{r-1}\gamma^{q^{r-1}-1},\gamma^{q^r-1}\right)=
\trdeg_{\F_p}\F_p\left
(a_0,\dots,a_{r-1}\gamma^{q^{r-1}-1},\gamma\right).
\end{equation}
On the other hand,
\begin{equation}
\label{E:t3}
\F_p(a_0,\dots,a_{r-1})\subset\F_p\left
(a_0,a_1\gamma^{q-1},\dots,a_{r-1}\gamma^{q^{r-1}-1},\gamma\right).
\end{equation}
Equations \eqref{E:t2} and \eqref{E:t3} yield
\eqref{E:t1}.
\end{proof}

\begin{definition}
\label{D:relmodtrdeg}
Let $K_0$ be any subfield of $K$. Then the relative
modular transcendence degree 
of $\phi$ over $K_0$ is the minimum transcendence
degree over $K_0$ of the 
compositum field of $K_0$ and the field of definition
of $\phi^{(\gamma)}$, the
minimum being taken over all $\gamma\in 
K^{\alg}\setminus\{0\}$. 
\end{definition}

Also, if for some non-constant $t\in A$, $\phi_t=\sum_{i=0}^r a_i\tau^i$ is monic, we can deduce that the relative modular transcendence degree of $\phi$ over $K_0$ can be 
defined as $$\trdeg_{K_0}K_0(a_{0},\dots,a_{r-1}),$$as an immediate corollary of Lemma~\ref{L:realmodtrdeg}.

\begin{theorem}
\label{T:T4mw}
Let $K$ be a countable field of characteristic $p$. Let $U$ be a coherent good set of valuations on $K$ and
let $F$ be the field of constants for $U$. Let
$\phi:A\rightarrow K\{\tau\}$ be a Drinfeld module of
positive relative modular transcendence degree over
$F$. Then $\phi(K)$ is a direct sum of a finite torsion
submodule and a free submodule of rank $\aleph_0$.
\end{theorem}

\begin{proof}
We first recall the definition of a \emph{tame}
module. The module $M$ is tame if every finite rank
submodule of $M$ is finitely generated. According to
Proposition $10$ from \cite{Poo}, in order to prove
Theorem~\ref{T:T4mw}, it suffices to show that
$\phi(K)$ is a tame module of rank $\aleph_0$.

We first prove the following lemma which will allow us to make certain reductions during the proof of Theorem~\ref{T:T4mw}.
\begin{lemma}
\label{L:finext}
Let $K'$ be a field extension of $K$. Assume that
$\phi(K')$ is a tame module of rank $\aleph_0$.
Then also $\phi(K)$ is a tame module of rank $\aleph_0$.
\end{lemma}

\begin{proof}[Proof of Lemma \ref{L:finext}.]
Let $K_0$ be the field of definition for $\phi$. By
Lemma~\ref{L:fieldofdeffingen}, $K_0$ is finitely
generated. Because $\phi$ has positive modular
transcendence degree, $\trdeg_{\F_p}K_0\ge 1$.
Thus, as proved in \cite{Wan}, $\phi(K_0)$ is a tame
module of rank $\aleph_0$. Thus $\phi(K)$ has rank $\aleph_0$ because both $\phi(K_0)$ and $\phi(K')$ have rank $\aleph_0$. Because $\phi(K')$ is tame, then every finite rank submodule of $\phi(K)\subset\phi(K')$ is finitely generated. Hence $\phi(K)$ is tame, as desired.  
\end{proof}

Let $t$ be a non-constant element of $A$. Let
$\phi_t=\sum_{i=0}^{r}a_i\tau^i$.

Let $\gamma\in K^{\alg}$ satisfy
$\gamma^{q^r-1}a_r=1$. Assume that
$\phi^{(\gamma)}\left(K(\gamma)\right)$ is a tame module
of rank $\aleph_0$. Because $\phi^{(\gamma)}$ is
isomorphic to $\phi$ over $K(\gamma)$, it follows that
also $\phi\left(K(\gamma)\right)$ is a tame module of rank
$\aleph_0$. Using Lemma~\ref{L:finext} for
$K'=K(\gamma)$, we obtain that $\phi(K)$ is a direct
sum of a finite torsion submodule and a free module of
rank $\aleph_0$. Thus, it suffices to prove Theorem
~\ref{T:T4mw} under the hypothesis that $\phi_t$ is
monic.

Because $F$ is the field of constants with respect to $U$, then $F$ is algebraically closed
in $K$. 

Let $S_0$ be the set of
places in $U$ where $\phi$ has bad reduction.
Because we supposed that $\phi_t$ is monic, Lemma~\ref{L:bad reduction} yields that $S_0$ is the set of
places from $U$ where not all of the
coefficients $a_0,\dots,a_{r-1}$ are integral. 

\begin{lemma}
\label{L:S_0notempty}
The set $S_0$ is not empty.
\end{lemma}

\begin{proof}[Proof of Lemma \ref{L:S_0notempty}.]
If $S_0$ is empty, then by Lemma \ref{L:Lconstants},
$a_i\in F$ for all $i$. Then by Lemma
~\ref{L:fieldef}, we derive that $\phi$ is defined
over $F^{\alg}\cap K=F$, which is a contradiction with
our assumption that $\phi$ has positive relative
modular transcendence degree over $F$.
\end{proof}

Because $S_0$ is not empty, we use Theorem
~\ref{T:T2mwg} $b)$ and conclude that for every
non-torsion $x\in K$, there exists $v\in U$ such
that 
\begin{equation}
\label{E:hkfvineq}
\hhat_{U,v}(x)>
q^{-r\left(2+(r^2+r)|S_0|\right)}d(v).
\end{equation}
Using inequality \eqref{E:important inequality}, we
conclude that
\begin{equation}
\label{E:hkfineq1}
\hhat_{U,v}(x)>
\frac{q^{-r\left(2+(r^2+r)|S_0|\right)}}{[K:F(x_1,\dots,x_n)]}=:c(\phi,K)=c>0.
\end{equation}
Because $\hhat_{U}(x)\ge\hhat_{U,v}(x)$ we
conclude that for every non-torsion $x\in K$,
\begin{equation}
\label{E:the_one}
\hhat_{U}(x)>c.
\end{equation}

On the other hand, Theorem \ref{T:T2mwg} $b)$ shows
that $\phi_{\tor}(K)$ is bounded. Moreover, if
$b(t)\in\F_q[t]$ is the least common multiple
of all polynomials in $t$ of degree at most
$(r^2+r)|S_0|$, then for every
$x\in\phi_{\tor}(K)$, $\phi_{b(t)}(x)=0$.

The last ingredient of our proof is the next lemma.

\begin{lemma}
\label{L:heightonmodules}
Let $R$ be a Dedekind domain and let $M$ be an
$R$-module. Assume there exists a function
$h:M\rightarrow\mathbb{R}_{\ge 0}$ satisfying the
following properties

(i) (triangle inequality) $h(x\pm y)\le h(x)+h(y)$,
for every $x,y\in M$.

(ii) if $x\in M_{\tor}$, then $h(x)=0$.

(iii) there exists $c>0$ such that for each $x\notin
M_{\tor}$, $h(x)>c$.

(iv) there exists $a\in R\setminus\{0\}$ such that $R/aR$ is finite and for all 
$x\in M$, $h(ax)\ge 4h(x)$.

If $M_{\tor}$ is finite, then $M$ is a tame
$R$-module.
\end{lemma}

We first show how Lemma~\ref{L:heightonmodules} yields Theorem~\ref{T:T4mw} and then prove Lemma~\ref{L:heightonmodules}.

Let $a\in A$ be an element such that $q^{\deg(\phi_a)}\ge 4$ (we will need this assumption because we will apply next Lemma~\ref{L:heightonmodules}).
Because of the finiteness of $\phi_{\tor}(K)$ and
because of the equation \eqref{E:the_one}, the Dedekind
domain $A$, $a\in A$, $\phi(K)$ and $\hhat_{U}$
satisfy the hypothesis of Lemma
~\ref{L:heightonmodules} (note that $A/aA$ is finite as shown in \cite{Poo}). We 
conclude that $\phi(K)$
is a tame module. Because $\phi(K)$ is countable, it has at most countable rank. On the other 
hand, as shown in the proof of Lemma~\ref{L:finext}, $\phi(K)$ has at least 
countable rank because $\phi$ has positive modular transcendence degree. Thus 
$\phi(K)$ has rank $\aleph_0$. Finally, Proposition $10$ of \cite{Poo} shows that a tame module of rank $\aleph_0$ is a direct sum of a finite torsion submodule and a free submodule of rank $\aleph_0$.

In order to finish the proof of our Theorem~\ref{T:T4mw} we need to prove Lemma~\ref{L:heightonmodules}.

\begin{proof}[Proof of Lemma~\ref{L:heightonmodules}.]
By the definition of a tame module, it suffices to
assume that $M$ is a finite rank $R$-module and
conclude that it is finitely generated.

Let $a\in R$ as in $(iv)$ of Lemma
\ref{L:heightonmodules}. By Lemma $3$ of \cite{Poo},
$M/aM$ is finite (here we use the assumption that
$M_{\tor}$ is finite). The following result is the key
to the proof of Lemma~\ref{L:heightonmodules}.

\begin{sublemma}
\label{S:finitely}
For every $D>0$, there exists finitely many $x\in M$
such that $h(x)\le D$.
\end{sublemma}

\begin{proof}[Proof of Sublemma \ref{S:finitely}.]

If we suppose Sublemma~\ref{S:finitely} is not true,
then we can define 
$$C=\inf\{D\mid\text{ there exists infinitely many
$x\in M$ such that 
$h(x)\le D$}\}.$$
Properties $(ii)$ and $(iii)$ and the finiteness of
$M_{\tor}$ yield $C\ge c>0$. By the definition of $C$,
it must be that there exists an infinite sequence of
elements $z_n$ of $M$ such that for every $n$, 
$$h(z_n)<\frac{3C}{2}.$$

Because $M/aM$ is finite, there exists a coset of $aM$
in $M$ containing infinitely many $z_n$ from the above
sequence.

But if $k_1\ne k_2$ and $z_{k_1}$ and $z_{k_2}$ are in
the same coset of $aM$ in $M$, then let 
$y\in M$ be such that $ay=z_{k_1}-z_{k_2}$. Using
properties $(iv)$ and $(i)$, we get
$$h(y)\le\frac{h(z_{k_1}-z_{k_2})}{4}\le\frac{h(z_{k_1})+h(z_{k_2})}{4}<\frac{3C
}{4}.$$
We can do this for any two elements of the sequence
that lie in the same 
coset of $aM$ in $M$. Because there are infinitely
many of them lying in 
the same coset, we can construct infinitely many
elements $z\in M$ such that $h(z)<\frac{3C}{4}$,
contradicting the minimality of $C$.
\end{proof}

From this point on, our proof of Lemma \ref{L:heightonmodules} follows the 
classical
descent argument in the 
Mordell-Weil theorem (see \cite{Ser}).

Take coset representatives $y_1,\dots,y_k$ for $aM$ in
$M$. Define then 
$$B=\max_{i\in\{1,\dots,k\}}h(y_i).$$
Consider the set $Z=\{x\in M\mid h(x)\le B\}$, which
is finite 
according to Sublemma~\ref{S:finitely}. Let 
$N$ be the finitely generated $R$-submodule of $M$
which is spanned by $Z$.

We claim that $M=N$. If we suppose this is not
the case, then by Sublemma~\ref{S:finitely} we can
pick $y\in M-N$ which minimizes $h(y)$. Because $N$
contains all the coset representatives of 
$aM$ in $M$, 
we can find $i\in\{1,\dots,k\}$ such that $y-y_i\in
aM$. Let $x\in M$ be 
such that
$y-y_i=ax$. Then $x\notin N$ because 
otherwise it would follow that $y\in N$ (we already know $y_i\in N$).
By our choice of $y$ and by properties $(iv)$ and
$(i)$, we have
$$h(y)\le
h(x)\le\frac{h(y-y_i)}{4}\le\frac{h(y)+h(y_i)}{4}
\le\frac{h(y)+B}{4}.$$
This means that $h(y)\le\frac{B}{3}<B$. This
contradicts the fact that 
$y\notin N$ because $N$ contains all the elements
$z\in M$ such that $h(z)\le B$. This contradiction
shows that indeed $M=N$ and so, $M$ is finitely 
generated. 
\end{proof}

Let $a\in A$ be an element such that $q^{\deg(\phi_a)}\ge 4$ (we will need this assumption because we will apply next Lemma~\ref{L:heightonmodules}).
Because of the finiteness of $\phi_{\tor}(K)$ and
because of the equation \eqref{E:the_one}, the Dedekind
domain $A$, $a\in A$, $\phi(K)$ and $\hhat_{U}$
satisfy the hypothesis of Lemma
~\ref{L:heightonmodules} (note that $A/aA$ is finite as shown in \cite{Poo}). We 
conclude that $\phi(K)$
is a tame module. Because $\phi(K)$ is countable, it has at most countable rank. On the other 
hand, as shown in the proof of Lemma~\ref{L:finext}, $\phi(K)$ has at least 
countable rank because $\phi$ has positive modular transcendence degree. Thus 
$\phi(K)$ has rank $\aleph_0$. Proposition $10$ of \cite{Poo} finishes the proof of Theorem~\ref{T:T4mw}.
\end{proof} 

The following result is an immediate corollary to Theorem~\ref{T:T4mw}.
\begin{theorem}
\label{T:T4mwfg}
Let $F$ be a countable field of characteristic $p$ and
let $K$ be a finitely 
generated field over $F$. Let $\phi:A\rightarrow
K\{\tau\}$ be a Drinfeld module 
of positive relative modular transcendence degree over
$F$. Then $\phi(K)$ is a direct sum
of a finite torsion 
submodule and a free submodule of rank $\aleph_0$.
\end{theorem}

\begin{proof}
The coherent good set $U$ of valuations on $K$ (from the statement of Theorem~\ref{T:T4mw}) is the set $M_{K/F}$ constructed in Section $3$.
\end{proof}

The following result gives a structure theorem for Drinfeld modules which are defined over the field of constants (with respect to some coherent good set of valuations).
\begin{theorem}
\label{T:T4_0}
Let $F$ be a countable, algebraically closed field of characteristic $p$ and let $K$ be a finitely 
generated extension of $F$ of positive transcendence degree over $F$. If $\phi:A\rightarrow F\{\tau\}$ is a Drinfeld 
module, then $\phi(K)$ is the direct sum of $\phi(F)$ and a free 
submodule of rank $\aleph_0$.
\end{theorem}

\begin{proof}
Let $t$ be a non-constant element of $A$. Because $\phi$ is defined over $F$ and $F$ is algebraically closed, we 
can
find $\gamma\in F$ such that $\phi^{(\gamma)}_t$ is
monic. Because $\phi$ and $\phi^{(\gamma)}$ are isomorphic over $F$, it 
suffices to prove Theorem~\ref{T:T4_0} for $\phi^{(\gamma)}$. Thus we assume 
from now on that $\phi_t$ is monic.

We will show next that the module
$\phi(K)/\phi(F)$
is tame.

Let $\{x_1,\dots,x_n\}$ be a transcendence basis for
$K/F$. We construct the good set of
valuations
$M_{K/F}$ with
respect to $\{x_1,\dots,x_n\}$, as described in
Section $3$. Then we construct the local and global
heights associated to $\phi$.

\begin{lemma}
\label{L:const}
For every $x\in F$, $\hhat_{K/F}(x)=0$.
\end{lemma}

\begin{proof}[Proof of Lemma \ref{L:const}.]
For every $x\in F$ and for every $a\in A$, because $\phi$ is defined over 
$F$, $\phi_a(x)\in F$. Hence $v(\phi_a(x))=0$ and so, for every $v\in 
M_{K/F}$, $\hhat_{K/F,v}(x)=0$.
\end{proof}

We define $\Hhat:\phi(K)/\phi(F)\rightarrow
\mathbb{R}_{\ge 0}$ by
$$\Hhat\left(x+\phi(F)\right)=\hhat_{K/F}(x)$$ for every
$x\in K$. We will prove in the next lemma that this newly defined function is indeed well-defined.

\begin{lemma}
\label{L:well-defined}
The function $\Hhat$ is well-defined.
\end{lemma}

\begin{proof}[Proof of Lemma \ref{L:well-defined}.]
To show that $\Hhat$ is well-defined, it suffices to
show that for every $x,y\in K$, if
$x-y=z\in F$, then
$\hhat_{K/F}(x)=\hhat_{K/F}(y)$.

Using the triangle inequality and using
$\hhat_{K/F}(z)=0$
(see Lemma~\ref{L:const}), we get
\begin{equation}
\label{E:hxy}
\hhat_{K/F}(x)\le\hhat_{K/F}(y)+\hhat_{K/F}(z)=
\hhat_{K/F}(y).
\end{equation}
Similarly, using this time
$\hhat_{K/F}(-z)=0$
(also $-z\in F$), we get
\begin{equation}
\label{E:hyx}
\hhat_{K/F}(y)\le\hhat_{K/F}(x)+\hhat_{K/F}
(-z)=\hhat_{K/F}(x).
\end{equation}
Inequalities \eqref{E:hxy} and \eqref{E:hyx} show that
$\hhat_{K/F}(x)=\hhat_{K/F}(y)$,
as desired.
\end{proof}

For each $x\in K$, we denote by $\xbar$ its image in $\phi(K)/\phi(F)$.
\begin{lemma}
\label{L:Hheight}
The function $\Hhat$ satisfies the properties:

(i) $\Hhat\left(\overline{x+y}\right)\le\Hhat\left(\xbar\right)+\Hhat\left(\ybar\right)$, for all $x,y\in K$.

(ii) $\Hhat\left(\phi_a(\xbar)
\right)=\deg(\phi_a)\cdot\Hhat\left(\xbar
\right)$, for all $x\in K$ and all
$a\in A\setminus\{0\}$.

(iii) $\Hhat\left(\xbar
\right)\ge
\frac{1}{\left[K:F(x_1,\dots,x_n)\right]}$,
for all $x\notin F$.
\end{lemma}

\begin{proof}[Proof of Lemma \ref{L:Hheight}.]
Properties (i) and (ii) follow immediately from the
definition of $\Hhat$ and the fact that $\phi$ is defined over $F$ and
$\hhat_{K/F}$
satisfies the triangle inequality and
$\hhat_{K/F}(\phi_a(x))=\deg(\phi_a)\cdot
\hhat_{K/F}(x)$,
for all $x\in K$ and all $a\in
A\setminus\{0\}$.

Using the result of Theorem \ref{T:T2mwg} part $a)$,
we conclude that if $x\notin F$,
there exists $v\in
M_{K/F}$ such
that 
\begin{equation}
\label{E:Hheight1}
\hhat_{K/F,v}(x)\ge
d(v).
\end{equation}
Using inequality \eqref{E:important inequality} in
\eqref{E:Hheight1}, we get
$\hhat_{K/F,v}(x)\ge
\frac{1}{\left[K:F(x_1,\dots,x_n)\right]}$.

Because
$\hhat_{K/F}(x)\ge\hhat_{K/F,v}(x)$,
we conclude that
$$\hhat_{K/F}(x)\ge
\frac{1}{\left[K:F(x_1,\dots,x_n)\right]}.
$$ 
\end{proof}

Now we can finish the proof of Theorem \ref{T:T4_0}.
The rank of
$\phi(K)/\phi(F)$
is at most $\aleph_0$ because $K$
is countable ($F$ is countable and $K$ is a finitely generated extension of 
$F$). We know that
$\phi(K)/\phi(F)$
is torsion-free (if $\phi_a(x)\in F$ for some $a\in A\setminus\{0\}$, 
then $x\in F$, because $\phi_a\in F\{\tau\}$). Because $\Hhat$ satisfies the properties
(i)-(iii) from Lemma~\ref{L:Hheight}, Lemma
~\ref{L:heightonmodules} yields that
$\phi(K)/\phi(F)$
is tame. 

\begin{lemma}
\label{L:rang}
The rank of $\phi(K)/\phi(F)$ is $\aleph_0$.
\end{lemma}

\begin{proof}[Proof of Lemma \ref{L:rang}.]
We need to show only that the rank of the above module is at least $\aleph_0$. 
Assume the rank is finite and we will derive a contradiction. 

Let $y_1,\dots,y_g\in K$ be the generators of 
$\left(\phi(K)/\phi(F)\right)\tensor_A\Frac(A)$ as a 
$\Frac(A)$-vector space. Let $v\in M_{K/F}$ be a place different 
from the finitely many places from $M_{K/F}$ where 
$y_1,\dots,y_g$ have poles. Let $x\in K$ be an element which has a pole 
at $v$. Then for every $a\in A\setminus\{0\}$, $\phi_a(x)$ has a pole at $v$. On 
the other hand, for every $a\in A$ and every $i\in\{1,\dots,g\}$, $\phi_a(y_i)$ 
is integral at $v$. Thus the equation 
$$\phi_a(x)=z+\sum_{i=1}^g\phi_{a_i}(y_i)$$
has no solutions $a,a_1,\dots,a_g\in A$ and $z\in F$ with $a\ne 0$. This 
provides a contradiction to our assumption that $y_1,\dots,y_g$ are generators 
for $\left(\phi(K)/\phi(F)\right)\tensor_A\Frac(A)$ as a 
$\Frac(A)$-vector space.
\end{proof}

Hence the rank of $\phi(K)/\phi(F)$ is $\aleph_0$.
Because $\phi(K)/\phi(F)$ is tame, Proposition $10$ of \cite{Poo} yields that
$\phi(K)/\phi(F)$
is a direct sum of its torsion submodule and a free
submodule of rank $\aleph_0$. As explained before,
$\phi(K)/\phi(F)$
is torsion-free. Hence
$\phi(K)/\phi(F)$
is free of rank $\aleph_0$. We have the exact
sequence:
$$0\rightarrow\phi(F)\rightarrow\phi(K)
\rightarrow\phi(K)/\phi(F)\rightarrow
0.$$
Because
$\phi(K)/\phi(F)$
is free, the above exact sequence splits. Thus,
$\phi(K)$ is a direct sum of
$\phi(F)$ and a free submodule of
rank $\aleph_0$. 
\end{proof}

The following result is an immediate corollary of Theorem \ref{T:T4_0}.
\begin{theorem}
\label{T:T4last}
Let $K$ be a finitely generated field of positive
transcendence degree over $\F_p$. If
$\phi:A\rightarrow K\{\tau\}$ is a Drinfeld module
defined over a finite subfield of $K$, then 
$\phi(\F_p^{\alg}K)$ is a direct sum of an
infinite torsion submodule 
(which is $\F_p^{\alg}$, the entire torsion
submodule of $\phi$) and a free submodule of rank
$\aleph_0$.
\end{theorem}

\section{Drinfeld modules over the perfect closure of a field}

In this section we will prove a similar result as Theorem~\ref{T:T4mw} valid for the perfect closure of the field $K$ (as always, $\phi:A\rightarrow K\{\tau\}$). Even though the result is an extension to Theorem~\ref{T:T4mw} and the general idea of its proof is similar with the one from Theorem~\ref{T:T4mw}, it makes more sense to be presented in a separate section. One reason is that it requires more refined height inequalities for Drinfeld modules as the ones proved so far. Also, the results of this section should be seen as an analogue of the author's results from \cite{ellper} (see also Chapter $3$ of \cite{teza}). In \cite{ellper} we proved a Mordell-Weil type theorem for non-isotrivial elliptic curves over the perfect closure of a function field of a curve over a finite field.

The setting for this section is the following: $K$ is a field of characteristic $p$ and $U$ is a coherent good set of valuations on $K$. Let $K_0\subset K$ be the field of constants with respect to $U$.

Let $\phi:A\rightarrow K\{\tau\}$ be a Drinfeld module. We construct the global height $\hhat$ and the local heights $\hhat_v$ with respect to the valuations in $U$ and the Drinfeld module $\phi$.

Assume $\phi$ has positive relative modular transcendence degree over $K_0$. Our goal is to prove there exists a constant $C>0$ depending only on $\phi$ and $K$ such that for every non-torsion point $x\in K^{\per}$, $\hhat(x)\ge C$. Clearly, it suffices to prove our result for an extension $L$ of $K$ (as long as we can control the dependence of the constant $C$ on the field extension). Also, replacing $\phi$ by an isomorphic Drinfeld module does not affect the validity of our statement. Therefore, we may assume as before, that for some non-constant $t\in A$, $\phi_t$ is monic.

Let $\phi_t=\sum_{i=0}^r a_i\tau^i$ (with $a_r=1$). Let $S_0$ be the set of places $v\in U$ for which there exists $i\in\{0,\dots,r\}$ such that $v(a_i)<0$. By Lemma~\ref{L:bad reduction}, $S_0$ is the finite set of places $v\in U$ of bad reduction for $\phi$. Not all of the coefficients $a_i$ are constant, because this would imply $\phi$ is defined over $K_0^{\alg}\cap K=K_0$. Therefore $S_0$ is not empty.

Let $L$ be any finite purely inseparable extension of $K$. Let $U_L$ be the set of places of $L$ which lie above the places from $U$. We use the convention, as always, that each valuation function is normalized so that its range equals $\mathbb{Z}$. We let $S$ be the finite set of places $w\in U_L$ which lie above places $v\in S_0$. Then $|S|=|S_0|>0$ (above each place from $S_0$ lies an unique place from $S$ because $L/K$ is purely inseparable).

In this section we will use again the definitions of $M_v$, $P_v$ and $R_v(\alpha)$ for $\alpha\in P_v$. So, we recall that for all $v\in S$, $|P_v|\le r+1$ and  for each $\alpha\in P_v$, $|R_v(\alpha)|\le q^r$. 

As stated in Lemma~\ref{L:L5}, for every $x\in L$ and every $v\in S$, if $v(x)\le 0$, then either
\begin{equation}
\label{E:eq1}
(v(x),\ac_{\pi_v}(x))\in P_v\times R_v(v(x))
\end{equation}
or $\hhat_v(x)\ge\frac{-d(v)M_v}{q^r}$.

Fix $v\in S$. Let $v_0$ be the place of $S_0$ lying below $v$. We define 
\begin{equation}
\label{E:eq2}
T_v:=-\min_{0\le i\le r}\frac{v(a_i)}{q^i}. 
\end{equation}
Because $v\in S$, $T_v>0$. Moreover, $T_v=-[L:K]\cdot\min_{0\le i\le r}\frac{v_0(a_i)}{q^i}$ (we used that $L/K$ is purely inseparable and so, $e(v|v_0)=[L:K]$). Thus 
\begin{equation}
\label{E:eq7}
T_v\ge\frac{[L:K]}{q^r}. 
\end{equation}

Let $P'_v$ be the set of all $0<\alpha\le T_v$ such that 
\begin{equation}
\label{E:eq3}
\min_{0\le i\le r}\left(v(a_i)+q^i\alpha\right)=:\alpha_1\in P_v.
\end{equation}
Because the function $f(y)=\min_i \left(v(a_i)+q^iy\right)$ is piecewise strictly increasing, we conclude $|P'_v|\le |P_v|\le r+1$ (for each $\alpha_1\in P_v$, there exists at most one $\alpha\in P'_v$ such that \eqref{E:eq3} holds). 

For each $\alpha\in P'_v$, we let $i_1,\dots,i_l$ be all the indices $i_j$ such that $$v(a_{i_j})+q^{i_j}\alpha=\min_{0\le i\le r}\left(v(a_i)+q^i\alpha\right).$$ We let $R_v(\alpha)$ be the set of all $\beta$ such that 
\begin{equation}
\label{E:eq4}
\sum_{1\le j\le l}\ac_{\pi_v}(a_{i_j})\beta^{q^{i_j}}\in R_v(\alpha_1).
\end{equation}
Because $|R_v(\alpha_1)|\le q^r$, $|R_v(\alpha)|\le q^{2r}$.

Let $P''_v$ be the set of all $0<\alpha\le T_v$ such that $-\alpha$ is a slope of a segment in the Newton polygon for $\phi_t$. For each $\alpha\in P''_v$, we let $i_1,\dots,i_l$ be all the indices $i_j$ such that $v(a_{i_j})+q^{i_j}\alpha=\min_{0\le i\le r}\left(v(a_i)+q^i\alpha\right)$. We let $R_v(\alpha)$ be the set of all $\beta$ such that 
\begin{equation}
\label{E:eq5}
\sum_{1\le j\le l}\ac_{\pi_v}(a_{i_j})\beta^{q^{i_j}}=0.
\end{equation}
We note that it might be that $\alpha\in P'_v\cap P''_v$. In that case, $R_v(\alpha)$ contains all $\beta$ satisfying both \eqref{E:eq4} and \eqref{E:eq5}. Therefore $|R_v(\alpha)|\le q^{2r}+q^r<q^{2(r+1)}$. 

Let $Q_v:=P_v\cup P'_v\cup P''_v$. Then $|Q_v|\le |P'_v|+|P_v\cup P''_v|\le (r+1)+(r+1)=2(r+1)$ (the cardinality of $P_v\cup P''_v$ is at most $r+1$ because there are at most $r$ segments in the Newton polygon for $\phi_t$ and besides the negatives of the slopes of the segments in the Newton polygon of $\phi_t$, only the number $0$ might be contained in $P_v\cup P''_v$).

The following result should be seen as an extension of Lemma~\ref{L:L5}.
\begin{lemma}
\label{L:another dichotomy}
Let $v\in S$ and let $x\in L$. Assume $v(x)\le T_v$. If $(v(x),\ac_{\pi_v}(x))\notin Q_v\times R_v(v(x))$, then $\hhat_v(x)\ge\frac{-d(v)M_v}{q^{2r}}$.
\end{lemma}

\begin{proof}
There are two cases: $v(x)\le 0$ and $0<v(x)\le T_v$.

We analyze the first case: $v(x)\le 0$. Because $(v(x),\ac_{\pi_v}(x))\notin P_v\times R_v(v(x))$, Lemma~\ref{L:L5} yields $\hhat_v(x)\ge\frac{-d(v)M_v}{q^r}>\frac{-d(v)M_v}{q^{2r}}$.

Assume now that $0<v(x)\le T_v$. Because $(v(x),\ac_{\pi_v}(x))\notin P''_v\times R_v(v(x))$, $v(\phi_t(x))=\min_{0\le i\le r}v(a_ix^{q^i})$ (see the remark at the end of the proof of Lemma~\ref{L:L0}). Because $v(x)\le T_v$, $v(a_ix^{q^i})\le 0$, for some $0\le i\le r$ (see the definition of $T_v$ from \eqref{E:eq2}). Hence $v(\phi_t(x))\le 0$. Let $i_1,\dots,i_l\in\{0,\dots,r\}$ be all the indices $i_j$ such that $v(\phi_t(x))=v(a_{i_j}x^{q^{i_j}})$. Then $\ac_{\pi_v}(\phi_t(x))=\sum_{j}\ac_{\pi_v}(a_{i_j})\ac_{\pi_v}(x)^{q^{i_j}}$. Because $(v(x),\ac_{\pi_v}(x))\notin P'_v\times R_v(v(x))$, we conclude 
\begin{equation}
\label{E:eq6}
(v(\phi_t(x)),\ac_{\pi_v}(\phi_t(x))\notin P_v\times R_v(v(\phi_t(x))).
\end{equation}
Because $v(\phi_t(x))\le 0$, Lemma~\ref{L:L5} yields $\hhat_v(\phi_t(x))\ge\frac{-d(v)M_v}{q^r}$. Hence $\hhat_v(x)=\frac{\hhat_v(\phi_t(x))}{q^r}\ge\frac{-d(v)M_v}{q^{2r}}$, as desired.
\end{proof}

The proof of following result is similar with the proof of Theorem~\ref{T:T2mwg}.
\begin{lemma}
\label{L:key dichotomy}
Let $x\in L$. Then either there exists $v\in S$ such that $\hhat_v(x)\ge\frac{-d(v)M_v}{q^{4r(r+1)^2|S|+2r}}$, or there exists a polynomial $b\in\F_q[t]$ of degree at most $4(r+1)^2|S|$ such that for every $v\in S$, $v(\phi_b(x))>T_v$.
\end{lemma}

\begin{proof}
Let $v\in S$. We apply Lemma~\ref{L:L-1} with $N=T_v$, $I=Q_v$ and $R(\alpha)=R_v(\alpha)$ for every $\alpha\in Q_v$. Because $|Q_v|\le 2(r+1)$ and $|R_v(\alpha)|<q^{2(r+1)}$, for every $\alpha\in Q_v$, we conclude that the following is true.

\begin{Fact}
\label{F:fact1}
Let $v\in S$. Let $W$ be an $\F_q$-vector subspace of $L$ with the property that for all $w\in W$, $(v(w),\ac_{\pi_v}(w))\in Q_v\times R_v(v(w))$ whenever $v(w)\le T_v$.

Then the $\F_q$-codimension of $\{w\in W\mid v(w)>T_v\}$ in $W$ is bounded above by $4(r+1)^2$.
\end{Fact}

We apply Fact~\ref{F:fact1} for each $v\in S$ and we deduce the following two results.

\begin{Fact}
\label{F:fact2}
Let $W$ be an $\F_q$-subspace of $L$ such that $(v(w),\ac_{\pi_v}(w))\in Q_v\times R_v(v(w))$ whenever $v\in S$, $w\in W$ and $v(w)\le T_v$. Then the $\F_q$-codimension of
$$\{w\in W\mid v(w)>T_v\text{ for all }v\in S\}$$
in $W$ is bounded above by $4(r+1)^2|S|$.
\end{Fact}

\begin{Fact}
\label{F:fact3}
Let notation be as in Fact~\ref{F:fact2}. If moreover, $\dim_{\F_q}W>4(r+1)^2|S|$, then there exists a nonzero $w\in W$ such that $v(w)>T_v$ for all $v\in S$.
\end{Fact}

We are now ready to finish the proof of Lemma~\ref{L:key dichotomy}. 

Let $W=\Span_{\F_q}\left(\{x,\phi_t(x),\dots,\phi_{t^{4(r+1)^2|S|}}\}\right)$. Because $\dim_{\F_q}W=4(r+1)^2|S|+1$, Fact~\ref{F:fact3} yields the existence of a nonzero $w\in W$ such that either there exists $v\in S$ such that 
\begin{equation}
\label{E:eq1dichotomy}
v(w)\le T_v\text{ and }(v(w) ,\ac_{\pi_v}(w))\notin Q_v\times R_v(v(w)),
\end{equation}
or for all $v\in S$,
\begin{equation}
\label{E:eq2dichotomy}
v(w)>T_v.
\end{equation}
If \eqref{E:eq1dichotomy} holds, then by Lemma~\ref{L:another dichotomy}, $\hhat_v(w)\ge\frac{-d(v)M_v}{q^{2r}}$. Because $w\in W\setminus\{0\}$, there exists a nonzero polynomial $b\in\F_q[t]$ of degree at most $4(r+1)^2|S|$ such that $w=\phi_b(x)$. Thus 
$$\hhat_v(x)=\frac{\hhat_v(w)}{\deg(\phi_b)}\ge\frac{-d(v)M_v}{q^{4r(r+1)^2|S|+2r}}.$$
If \eqref{E:eq2dichotomy} holds, then $v(\phi_b(x))>T_v$ for every $v\in S$ with $\deg(b)\le 4(r+1)^2|S|$, where $w=\phi_b(x)$.
\end{proof}

The next theorem is the main step in order to prove the existence of a positive lower bound $C$ for the height of non-torsion points in $K^{\per}$ for a Drinfeld modules $\phi:A\rightarrow K\{\tau\}$ of positive relative modular transcendence degree. 

\begin{theorem}
\label{T:lehper}
Let $K$ be a field of chracteristic $p$ and let $U$ be a coherent good set of valuations on $K$. Let $K_0\subset K$ be the field of constants with respect to the valuations in $U$.  

Let $\phi:A\rightarrow K\{\tau\}$ be a Drinfeld module of positive relative modular transcendence degree over $K_0$. Let $t\in A$ be a non-constant element and assume $\phi_t=\sum_{i=0}^ra_i\tau^i$ is monic. Let $S_0$ be the set of the places in $U$ of bad reduction for $\phi$. Let $s=|S_0|$.

Then for every non-torsion point $x\in K^{\per}$, $\hhat(x)>\frac{\min_{v_0\in S_0}d(v_0)}{q^{4r(r+1)^2s+3r}}>0$. 
\end{theorem}

Before proving Theorem~\ref{T:lehper} we show how this theorem implies a more general result.
\begin{theorem}
\label{T:lehpergeneral}
Let $K$ be a field of chracteristic $p$ and let $U$ be a coherent good set of valuations on $K$. Let $K_0\subset K$ be the field of constants with respect to the valuations in $U$. Let $\phi:A\rightarrow K\{\tau\}$ be a Drinfeld module of positive relative modular transcendence degree over $K_0$. There exists a positive constant $C$ depending only on $\phi$ and $K$ such that for every non-torsion $x\in K^{\per}$, $\hhat(x)>C$.
\end{theorem}

\begin{proof}
Let $t\in A$ be a non-constant element minimizing $\deg(\phi_a)$ as $a$ ranges over non-constant elements of $A$. Let $\gamma\in K^{\alg}$ such that $\phi^{(\gamma)}_t$ is monic. Let $L=K(\gamma)$. If $\deg(\phi_t)=q^r$, then $[L:K]\le q^r-1$. As explained before, $\phi$ and $\phi^{(\gamma)}$ are isomorphic and $\hhat_{\phi}=\hhat_{\phi^{(\gamma)}}$.

Let $S$ be the set of all places of $L$ lying above places in $S_0$ (we use the notation from Theorem~\ref{T:lehper}). Then 
\begin{equation}
\label{E:extensie1}
|S|\le (q^r-1)|S_0|, 
\end{equation}
because $[L:K]\le q^r-1$. Also, for each $v\in S$, 
\begin{equation}
\label{E:extensie2}
d(v)=\frac{d(v_0)f(v|v_0)}{[L:K]}\ge\frac{d(v_0)}{[L:K]}.
\end{equation}
By Theorem~\ref{T:lehper} applied to $\phi^{(\gamma)}:A\rightarrow L\{\tau\}$, for every $x\in K^{\per}\subset L^{\per}$, 
\begin{equation}
\label{E:extensie3}
\hhat_{\phi}(x)=\hhat_{\phi^{(\gamma)}}(x)\ge\frac{\min_{v\in S}d(v)}{q^{4r(r+1)^2|S|+3r}}=:C>0.
\end{equation}
Using \eqref{E:extensie1} and \eqref{E:extensie2} in \eqref{E:extensie3}, we see that the positive constant $C$ is bounded below by another positive constant which depends only on $\phi$ and $K$, as desired.
\end{proof}

\begin{proof}[Proof of Theorem~\ref{T:lehper}.]
We first recall that because $\phi$ has positive relative modular transcendence degree over $K_0$, $s\ge 1$.

Let $x\in K^{\per}$ be a non-torsion point and let $L=K(x)$. Let $S$ be the set of places of $L$ which lie above places in $S_0$. Because $L\subset K^{\per}$, $|S|=s$.

Let $U_L$ be the set of places of $L$ which lie above the places from $U$ (as always, they are normalized so that the range of each valuation function is $\mathbb{Z}$). For each place $v\in U_L$, we denote by $v_0\in S_0$ the corresponding place which lies below $v$. Because $L/K$ is a purely inseparable extension, for each place $v\in U_L$, 
\begin{equation}
\label{E:eq8}
d(v)e(v|v_0)=\frac{f(v|v_0)d(v_0)e(v|v_0)}{[L:K]}=\frac{d(v_0)[L:K]}{[L:K]}=d(v_0).
\end{equation}
On the other hand, by its definition, $M_v=e(v|v_0)\min_{0\le i<r}\frac{v_0(a_i)}{q^r-q^i}$ and so, if $v\in S$, 
\begin{equation}
\label{E:eq9}
M_v<\frac{-e(v|v_0)}{q^r}.
\end{equation}
Using \eqref{E:eq8} in \eqref{E:eq9} we get $-d(v)M_v>\frac{d(v_0)}{q^r}$. 

If there exists $v\in S$ such that $\hhat_v(x)>\frac{d(v_0)}{q^{4r(r+1)^2s+3r}}$, then $$\hhat(x)\ge\hhat_v(x)>\min_{v_0\in S_0}\frac{d(v_0)}{q^{4r(r+1)^2s+3r}},$$ 
as desired. Therefore, assume from now on, that for each $v\in S$, $$\hhat_v(x)\le\frac{d(v_0)}{q^{4r(r+1)^2s+3r}}<\frac{-d(v)M_v}{q^{4r(r+1)^2s+2r}}.$$ 
Then by Lemma~\ref{L:key dichotomy}, there exists a nonzero polynomial $b\in\F_q[t]$ of degree at most $4(r+1)^2s$ such that 
\begin{equation}
\label{E:eq91}
v(\phi_b(x))>T_v
\end{equation}
for all $v\in S$. Because $b\ne 0$ and $x\notin\phi_{\tor}$, $y:=\phi_b(x)\ne 0$. Because $U$ is a coherent good set of valuations, $U_L$ is a good set of valuations on $L$ and so, because $y\ne 0$,
\begin{equation}
\label{E:eq10}
\sum_{v\in U_L}d(v)\cdot v(y)=0.
\end{equation}

By its definition, for each $v\in S$, $T_v=-e(v|v_0)\min_{0\le i\le r}\frac{v_0(a_i)}{q^i}\ge\frac{e(v|v_0)}{q^r}$. Hence, using \eqref{E:eq91} and \eqref{E:eq8}, we get 
\begin{equation}
\label{E:eq12}
\sum_{v\in S}d(v)\cdot v(y)>\sum_{v\in S}\frac{d(v)e(v|v_0)}{q^r}=\sum_{v_0\in S_0}\frac{d(v_0)}{q^r}.
\end{equation}

Using \eqref{E:eq12} in \eqref{E:eq10}, we conclude there exists a finite set $U(y)$ of places in $U_L\setminus S$ such that for each $v\in U(y)$, $v(y)<0$ and moreover
\begin{equation}
\label{E:eq11}
\sum_{v\in U(y)}d(v)\cdot v(y)<-\sum_{v_0\in S_0}\frac{d(v_0)}{q^r}.
\end{equation}

For each $v\in U(y)$, because $v\notin S$ and $v(y)<0$, Lemma~\ref{L:L1'} yields $$\hhat_v(y)=-d(v)\cdot v(y).$$ 
Using \eqref{E:eq11} we conclude
\begin{equation}
\label{E:eq13}
\sum_{v\in U(y)}\hhat_v(y)>\sum_{v_0\in S_0}\frac{d(v_0)}{q^r}.
\end{equation}

Inequality \eqref{E:eq13} yields 
\begin{equation}
\label{E:eq15}
\hhat(y)>\sum_{v_0\in S_0}\frac{d(v_0)}{q^r}.
\end{equation}
Because $y=\phi_b(x)$ and the degree of $b$ as a polynomial in $t$ is at most $4(r+1)^2s$, we get 
$$\hhat(x)\ge\frac{\hhat(y)}{q^{4r(r+1)^2s}}>\sum_{v_0\in S_0}\frac{d(v_0)}{q^{4r(r+1)^2s+r}}>\min_{v_0\in S_0}\frac{d(v_0)}{q^{4r(r+1)^2s+3r}},$$
as desired.
\end{proof}

\begin{remark}
\label{R:lehper is sharp}
Theorem~\ref{T:lehper} is sharp in the sense that if we assume $\phi:A\rightarrow K_0\{\tau\}$ and we keep the rest of the assumptions from Theorem~\ref{T:lehper}, then the conclusion of Theorem~\ref{T:lehper} fails. Indeed, let $x\in K\setminus K_0$. By Lemma~\ref{L:L1'}, $\hhat(x)=-\sum_{v\in U}d(v)\cdot\min\{0,v(x)\}$. Because $x\notin K_0$, $\hhat(x)>0$. Moreover, for each $n\ge 1$, $x^{1/p^n}\in K^{\per}$ and an easy computation, using again Lemma~\ref{L:L1'} shows $\hhat(x^{1/p^n})=\frac{\hhat(x)}{p^n}$. Therefore, as $n$ goes to infinity, the height of $x^{1/p^n}$ goes strictly decreasing to $0$. Hence there is no uniform positive lower bound for the height of a non-torsion point in $K^{\per}$.
\end{remark}

\begin{corollary}
\label{T:lehperfunc}
Let $K$ be a finitely generated field of characteristic $p$. Let $\phi:A\rightarrow K\{\tau\}$ be a Drinfeld module of positive modular transcendence degree. There exists a constant $C>0$ depending only on $\phi$ and $K$ such that for every non-torsion point $x\in K^{\per}$, $\hhat(x)\ge C$.
\end{corollary}

\begin{proof}
Let $V$ be a projective normal variety defined over a finite field, whose function field is $K$. We construct the coherent good set $U$ of valuations on $K$ associated to the irreducible divisors of $V$ (in Section $3$ we presented a completely algebraic construction of $U$). The field of constants with respect to $U$ is the maximal finite subfield of $K$. Because $\phi$ has positive modular transcendence degree we can apply Theorem~\ref{T:lehpergeneral} and get the existence of the constant $C$ in Corollary~\ref{T:lehperfunc}.
\end{proof}

Using Theorem~\ref{T:lehper} we prove the following Mordell-Weil type theorem.
\begin{theorem}
\label{T:drpermw}
Let $K$ be a countable field of chracteristic $p$ and let $U$ be a coherent good set of valuations on $K$. Let $K_0\subset K$ be the field of constants with respect to the valuations from $U$. Let $\phi:A\rightarrow K\{\tau\}$ be a Drinfeld module of positive relative modular transcendence degree over $K_0$. Then $\phi(K^{\per})$ is the direct sum of a finite torsion submodule with a free submodule of rank $\aleph_0$.
\end{theorem}

\begin{proof}
We will prove $\phi(K^{\per})$ has rank $\aleph_0$ and is a tame module. According to Proposition $10$ of \cite{Poo}, these two properties yield our conclusion.

As proved in Lemma~\ref{L:finext}, it suffices to prove our theorem after replacing $K$ by a finite extension. Therefore, we assume from now on that there exists a non-constant $t\in A$ such that $\phi_t$ is monic. Let $q^r$ be the degree of $\phi_t$.

We know that $\phi(K^{\per})$ has rank at most $\aleph_0$ because $K^{\per}$ is countable, as $K$ is countable. On the other hand, $\phi(K^{\per})$ has at least rank $\aleph_0$ because $\phi$ has positive modular transcendence degree (see the proof of Lemma~\ref{L:finext}). 

In order to show $\phi(K^{\per})$ is tame, we use Lemma~\ref{L:heightonmodules}. Thus we need to show that $\phi_{\tor}(K^{\per})$ is finite. The other conditions of the above mentioned lemma are already satisfied by the global height function associated to $\phi$ (see Theorem~\ref{T:lehper}) and by any element $a\in A$ of degree at least $2$ (so that $\deg(\phi_a)=q^{\deg(a)}\ge 4$). Because $K^{\per}=\bigcup_{n\ge 1}K^{1/p^n}$, it suffices to show $\phi_{\tor}(K^{1/p^n})$ is uniformly bounded.

Let $s\ge 1$ be the number of places in $U$ of bad reduction for $\phi$. Let $U_n$ be the good set of places on $K^{1/p^n}$, which lie above places in $U$. There exists exactly one place in $U_n$ lying above each place in $U$ because $K^{1/p^n}/K$ is a purely inseparable extension. Thus, for each $n\ge 1$, there are $s$ places of bad reduction for $\phi$ in $U_n$. By Theorem~\ref{T:T2mwg}, the size of $\phi_{\tor}(K^{1/p^n})$ is bounded above in terms of $q$, $r$ and $s$, independently of $n$. Hence $\phi_{\tor}(K^{\per})$ is finite. As explained in the previous paragraph, Lemma~\ref{L:heightonmodules} concludes the proof of our theorem.
\end{proof}

Just as Corollary~\ref{T:lehperfunc} followed from Theorem~\ref{T:lehper}, in the same way we can deduce the following result from Theorem~\ref{T:drpermw}.
\begin{corollary}
\label{T:drpermwfunc}
Let $K$ be a finitely generated field of charateristic $p$ and let $\phi:A\rightarrow K\{\tau\}$ be a Drinfeld module of positive modular transcendence degree. Then $\phi(K^{\per})$ is a direct sum of a finite torsion submodule with a free submodule of rank $\aleph_0$.
\end{corollary}

\begin{remark}
\label{T:drpermwfunc is sharp}
Corollary~\ref{T:drpermwfunc} is sharp in the sense that we cannot drop the hypothesis that $\phi$ has positive modular transcendence degree. For example, let $\phi:\F_q[t]\rightarrow\F_q(t)\{\tau\}$ be given by $\phi_t=\tau$. Then we can check immediately that $\phi\left(\F_q(t)^{\per}\right)$ is a direct sum of a finite torsion submodule with a free $\F_q[t,t^{-1}]$-submodule of rank $\aleph_0$.
\end{remark}

\end{document}